\documentclass[a4paper]{amsart}

\usepackage[utf8]{inputenc}
\usepackage[T1]{fontenc}
\usepackage{lmodern, enumerate}
\usepackage{amssymb,amsxtra}
\usepackage{amscd}
\usepackage[all]{xy}
\usepackage{xcolor}
\usepackage{nicefrac,mathtools}
\usepackage{microtype}
\usepackage[pdftitle={Universal and exotic generalized fixed-point algebras},
 pdfauthor={Alcides Buss, Siegfried Echterhoff},
 pdfsubject={Mathematics}
]{hyperref}
\usepackage[lite]{amsrefs}

\numberwithin{equation}{section}
\theoremstyle{plain}
\newtheorem{theorem}[equation]{Theorem}
\newtheorem{lemma}[equation]{Lemma}
\newtheorem{proposition}[equation]{Proposition}
\newtheorem{corollary}[equation]{Corollary}
\theoremstyle{definition}
\newtheorem{definition}[equation]{Definition}
\theoremstyle{remark}
\newtheorem{remark}[equation]{Remark}
\newtheorem{example}[equation]{Example}

\DeclareMathOperator{\Aut}{Aut}

\DeclareMathOperator{\Ind}{Ind}
\DeclareMathOperator{\Prim}{Prim}

\DeclareMathOperator{\id}{\mathrm{id}}

\DeclareMathOperator{\Inf}{Inf}

\DeclareMathOperator{\Res}{\mathrm{Res}}

\DeclareMathOperator{\rank}{\mathrm{rank}}

\newcommand*{\R}{\mathbb R}
\newcommand*{\T}{\mathbb T}
\newcommand*{\C}{\mathbb C}
\newcommand*{\Z}{\mathbb Z}
\newcommand*{\N}{\mathbb N}
\renewcommand*{\L}{\mathcal L}
\newcommand*{\K}{\mathcal K}
\renewcommand*{\H}{\mathcal H}
\newcommand{\lk}{\langle}
\newcommand{\rk}{\rangle}
\newcommand{\lt}{\mathrm{lt}}
\newcommand{\rt}{\mathrm{rt}}

\newcommand*{\M}{\mathcal M}
\newcommand*{\X}{\mathcal X}

\newcommand{\dach}{{\widehat{ \ }}}

\newcommand*{\Ad}{\textup{Ad}}

\newcommand*{\U}{\mathcal U}
\newcommand*{\into}{\hookrightarrow}

\newcommand{\om}{\omega}

\newcommand{\car}{\curvearrowright}

\newcommand{\tr}{\mathrm{tr}}
\newcommand{\Q}{\mathbb{Q}}

\begin{document}
\title[Crossed products and twisted group algebras]{Simple subquotients of crossed products by abelian groups and twisted group algebras}

\author{Siegfried Echterhoff}
\email{echters@uni-muenster.de}
\address{Mathematisches Institut\\
Universit\"at M\"un\-ster\\
 Einsteinstr.\ 62\\
 48149 M\"unster\\
 Germany}

\begin{abstract}
Motivated by work of Poguntke we study the question under what conditions simple subquotients of crossed products 
$A\rtimes_{\alpha}G$  by (twisted) actions of abelian groups $G$  are isomorphic 
to simple twisted group algebras of abelian groups. As a consequence, we recover a theorem of Poguntke's  saying 
that the simple subquotients of group $C^*$-algebras of connected groups are either stably isomorphic to  $\C$ or 
they are stably isomorphic to simple non-commutative tori.
\end{abstract}



\thanks{Funded by the Deutsche Forschungsgemeinschaft (DFG, German Research Foundation) under Germany's Excellence Strategy EXC 2044/2 –390685587, Mathematics Münster: Dynamics–Geometry–Structure.}

\maketitle

\section{Introduction}\label{sec-intro}
One important motivation to study the structure of simple $C^*$-algebras $A$ comes from the fact that they are commonly considered as the building blocks of 
general $C^*$-algebras. Indeed, under the mild condition that points in the primitive ideal space $\Prim(A)$ are locally closed, 
each primitive ideal $P\subseteq A$ 
determines a unique simple subquotient $A_P:=\ker{\overline{\{P\}}}/P$ of $A$, where for $E\subseteq \Prim(A)$ we write  $\ker{E}:=\cap_{Q\in E}Q$. 
Then $A$ could be understood as an algebra of sections of a (usually very complicated) 
field of $C^*$-algebras over $\Prim(A)$ with fibres $A_P$. Of course, in general we are far away from a general understanding of this bundle structure even in the 
case where $\Prim(A)$ is Hausdorff, in which case we always obtain a continuous bundle of $C^*$-algebras. Nevertheless, given a $C^*$-algebra $A$, as a first step
it is  important to understand the structure of the fibres $A_P$, before we can even start to look at the global structure of the bundle.
 
The best behaved situation appears when $A$ is a $C^*$-algebra of type I, which just means that all fibres $A_P$ of $A$ are algebras $\K(\H)$ of compact operators on some Hilbert space $\H$.
Beyond that, the class of  simple non-commutative tori $A_\Theta$ is one  of the best studied  classes of simple $C^*$-algebras. Recall that, given  a skew symmetric matrix $\Theta\in M_n(\R)$, the $n$-dimensional 
non-commutative torus $A_\Theta$ is defined as the universal $C^*$-algebra generated by unitaries $u_1,\ldots, u_n$ subject to the relations 
$$u_i u_j=\exp(2\pi i \theta_{ij}) u_ju_i.$$
These algebras can also be realized as a twisted group algebra $A_\om=C^*(\Z^n,\om)$ by some group $2$-cocycle $\om\in Z^2(\Z^n,\T)$.
Since the isomorphism class of $C^*(\Z^n,\om)$ only depends on the cohomology class $[\om]\in H^2(\Z^n,\T)$, and since every 
cocycle $\om\in Z^2(\Z^n,\T)$ is cohomologous to one of the form $\om_\Theta(k,l)=\exp(2\pi i \lk \Theta k,l\rk), \quad k,l\in \Z^n,$
with  $\Theta$ as above,  this gives the same class of $C^*$-algebras. 
A non-commutative torus $A_\om$ is simple if and only if $\om$ is {\em totally skew} in the sense that the symmetry group 
 $$S_\om:=\{k\in \Z^n: \om(k,l)=\om(l,k)\;\forall l\in \Z^n\}$$
 is trivial. We have  the following remarkable theorem
of  Poguntke's  \cite{Pogi}*{Theorem 1 \& Theorem  2}, which  relies  on 
 some deep results due to Dixmier and Pukanszky on the representation theory of real algebraic groups and of connected groups:

\begin{theorem}[Poguntke]\label{thm-Pogi}
Suppose that $G$ is a  locally compact group.
 If $G$ is connected or compactly generated two-step nilpotent, then  every one-point set $\{P\}\in \Prim(C^*(G))$ is locally closed (even closed in the nilpotent case)
and the simple subquotients $C^*(G)_P$ of $C^*(G)$ are 
either Morita equivalent to $\C$ or to a simple non-commutative torus $A_\Theta$.
\end{theorem}

This work was mainly motivated by an attempt to better understand Poguntke's arguments and to put them in a somewhat broader
perspective. To start with, assume $G$ is a  connected and simply connected Lie group with Lie algebra $\mathfrak g$. 
By an application of Ado's theorem, we  can embed $\mathfrak g$ as an ideal of some real algebraic Lie algebra $\mathfrak l$
such that $\mathfrak n:=[\mathfrak g,\mathfrak g]=[\mathfrak l,\mathfrak l]$  (see the discussion in  \cite{Puk}*{\S 1}).
Let $N=\exp(\mathfrak n)$ and $L=\exp(\mathfrak l)$ denote the simply connected and connected Lie groups  associated to $\mathfrak n$
and $\mathfrak l$, respectively.  Then $N\subseteq G\subseteq L$ and   $N=[G,G]=[L,L]$. Since $N$ and $L$ are locally isomorphic to real algebraic groups, it follows from Dixmier's \cite{Dix}*{\S 2.1., proposition on p.~425} that $C^*(N)$ and $C^*(L)$ are of type I. 
Moreover,  Pukanszky's   \cite{Puk}*{Lemma 1.1.3} shows that the orbit space $\widehat{N}/L$ under the conjugation action
is countably separated, which, by Glimm's  \cite{Glimm}*{Theorem} is equivalent to say that the action 
$L\car \widehat{N}\cong \Prim(C^*(N))$  is smooth in the sense of Definition \ref{def-smooth} below. 

 Using Green's theory of twisted actions 
and twisted crossed products, there exists a twisted action $(\beta,\tau):(L,N)\car C^*(N)$ such that the $C^*$-group algebras 
$C^*(L)$ and $C^*(G)$ can be written as a twisted crossed product
$$ C^*(L)\cong C^*(N)\rtimes_{(\beta,\tau)}L/N\quad\text{and}\quad C^*(G)\cong C^*(N)\rtimes_{(\alpha,\tau)}G/N,$$
where $\alpha$ denotes the restriction of $\beta$ to $G$.
Basically starting from this,  Poguntke's proof of \cite{Pogi}*{Theorem 2} is given  in a series of quite technical and nontrivial reduction steps  to show that 
every primitive ideal $P\in \Prim(C^*(G))$ `lives' on a subquotient $B=J/I$ of $C^*(G)$ that is Morita equivalent 
to a twisted group algebra  $C^*(H,\om)$ for some compactly generated locally compact abelian group $H$ and some 
Borel $2$-cocycle $\om\in Z^2(H,\T)$. The statement on simple suquotients of $C^*(G)$ then follows from 
\cite{Pogi}*{Theorem 1}, which covers the nilpotent case of Theorem \ref{thm-Pogi} and implies that simple quotients of $C^*(H,\om)$ are either elementary (i.e., compact operators on some Hilbert  space) or  Morita equivalent to non-commutative tori.

In this paper we want to approach the above arguments by studying the general situation of crossed products by actions
$\alpha:G\car A$ of locally compact abelian groups $G$ on $C^*$-algebras $A$ such that 
$A\rtimes_\alpha G$ decomposes (in a certain sense) into a family of subquotients that are Morita equivalent to twisted group algebras 
of abelian groups. Recall that by \cite{DixC}*{Section 3.2} we may canonically identify the primitive ideal space of an ideal $I$ 
in a $C^*$-algebra $A$ with the open  subset $U_I:=\{P\in \Prim(A): I\not\subseteq P\}$ of $\Prim(A)$, and the primitive ideal space 
of the quotient $\Prim(A/I)$ with the closed subset $A_I:=\{P\in \Prim(A): I\subseteq P\}$. It then follows that if $J\subseteq I\subseteq A$ are
two ideals, then $\Prim(I/J)$ can be identified with the locally closed\footnote{meaning that $Y_{I/J}$ it is  relatively open in its closure}   set
$$Y_{I/J}:=\{P\in \Prim(A): I\not\subseteq P \;\text{and}\; J\subseteq P\}\subseteq \Prim(A).$$
Using this, the notion of `decomposition' we use here is the following:

\begin{definition}\label{def-decom}
Let $A$ be a $C^*$-algebra. For an index set $\Lambda$ let $(A_\lambda=I_\lambda/J_\lambda)_{\lambda\in \Lambda}$, be a family 
of subquotients of $A$. We say that $A$ {\em decomposes into the family of subquotients} $(A_\lambda)_{\lambda\in \Lambda}$ if
$$\Prim(A)=\dot{\bigcup}_{\lambda}\Prim(A_\lambda)$$
is the disjoint union of the locally closed subsets $\Prim(A_\lambda)$ of $\Prim(A)$.
\end{definition}

Let us  also recall the notion  of a smooth action. Recall that if $\alpha:L\car A$ is an action, then $L$ act on $\Prim(A)$ via 
$(g,P)\mapsto \alpha_g(P)$.

\begin{definition}\label{def-smooth} Suppose that $\alpha:L\car A$ is an action of a locally compact  group $L$
on a $C^*$-algebra  $A$. We say the action is {\em smooth}, if the following hold:
\begin{enumerate}
\item the orbit space $\Prim(A)/L$ is almost Hausdorff\footnote{Recall that a topological space $Y$ is called almost Hausdorff  if every closed subset $Z\subseteq Y$ contains a dense relatively open Hausdorff subset $U\subseteq Z$.} or $A$ is separable and $L$ is second countable, and
\item for all $P\in \Prim(A)$ the  orbit $L(P)=\{\alpha_g(L): g\in L\}$ is locally closed in $\Prim(A)$ and the map
$$L/L_P\to L(P); \; gL_P\mapsto \alpha_g(P)$$
is a homeomorphism, where $L_P=\{g\in L: \alpha_g(P)=P\}$ denotes the stabilizer of $P$.
\end{enumerate}
\end{definition}

Recall from   \cite{DixC}*{Chapter 4} that for every type I $C^*$-algebra $A$ we have $\widehat{A}\cong \Prim(A)$ is always almost Hausdorff.   By \cite{Glimm}*{Theorem} we know that if $A$ is separable and of type I, and $L$ is second countable, then  $\alpha:L\car A$ is 
smooth if and only if   $\widehat{A}/L$ is  a T$_0$-space, which always holds if $\widehat{A}/L$ is almost Hausdorff.
Moreover, it follows from \cite{Ungermann}*{Theorem 27} that actions of compact groups on $C^*$-algebras with almost Hausdorff 
primitive ideal spaces, in particular on $C^*$-algebras of type I,  are always smooth.

Our main result,  motivated by Poguntke's work,   uses the 
dual action $\widehat\beta:\widehat{L}\car A\rtimes_\beta L$ of the dual group $\widehat{L}$ 
on a crossed product $A\rtimes_\beta L$ of an abelian group $L$.

\begin{theorem}\label{thm-main-intro}
Suppose that $\beta:L\car A$ is an action of the abelian group $L$ on a $C^*$-algebra $A$ such that the following assumptions hold
\begin{enumerate}
\item[(i)] $A\rtimes_\beta L$ is type I, and 
\item[(ii)] the dual action $\widehat\beta:\widehat{L}\car A\rtimes_\beta L$ is smooth.
\end{enumerate}
Suppose further that the abelian group $G$ acts on $A$ via an action $\alpha=\beta\circ \varphi:G\car A$ for some continuous homomorphism
$\varphi:G\to L$. 
Then  there exists a decomposition of $A\rtimes_\alpha G$ by subquotients $(B_\lambda)_{\lambda \in \Lambda}$
with  index set $\Lambda =\Prim(A)/L$, such that each $B_\lambda$
is Morita equivalent to some twisted group algebra $C^*(H_\lambda,\om_\lambda)$ of some abelian locally compact group  $H_\lambda$
and some $2$-cocycle $\om_\lambda\in Z^2(H_\lambda, \T)$. 
\end{theorem}

By the structure theory of twisted group algebras of abelian groups due to Baggett and Kleppner (\cite{BagKlep}) it follows that
the primitive ideal spaces  $\Prim(C^*(H_\lambda,\om_\lambda))$ are always Hausdorff, and each simple quotient of 
$C^*(H_\lambda,\om_\lambda)$ is itself isomorphic to a twisted group algebra $C^*(H_\lambda/S_\lambda, \tilde\om_{\lambda})$
of some quotient group $H_\lambda/S_\lambda$ of $H_\lambda$. If $H_\lambda/S_\lambda$ is compactly generated,  Poguntke's \cite{Pogi}*{Theorem 1} (or the results in Section \ref{sec-twisted-group-algebras} below) show that 
$C^*(H_\lambda/S_\lambda, \tilde\om_{\lambda})$ is either Morita equivalent to $\C$ or to a simple non-commutative torus.

\begin{corollary}\label{cor-main-intro}
Suppose that $\beta:L\car A$ and $\alpha=\beta\circ \varphi:G\car A$ are  as in Theorem \ref{thm-main-intro}.
Then every primitive ideal $P\in \Prim(A\rtimes_\alpha G)$ is locally closed
and the corresponding subquotient $B_P$ of $A\rtimes_\alpha G$ is Morita equivalent to a simple twisted group algebra 
$C^*(H_P,\om_P)$ of some abelian group $H_P$. Moreover, if in addition $G$ is compactly generated and $L$ is a Lie group,
then the $B_P$ are either Morita equivalent to $\C$ or to some non-commutative torus $C^*(\Z^l,\om_\Theta)$ of dimension $l\geq 2$.
\end{corollary}

We shall see  that  conditions (i) and (ii) of Theorem \ref{thm-main-intro} always hold  if
$A$ and $A\rtimes_\beta L$ are both type I and $\beta:L\car A$ or  the dual action $\widehat\beta:\widehat{L}\car A\rtimes_\beta L$ is smooth. In the setting of twisted 
actions, this is precisely what appears to be true for the twisted action 
$(\beta,\tau):L/N\car C^*(N)$ in the above discussion.
Using the fact that every twisted action is Morita equivalent to an ordinary action, and that everything 
in the statement of Theorem \ref{thm-main-intro} is preserved under passing to Morita equivalent 
(twisted) actions, we see that  Theorem \ref{thm-Pogi} can be obtained as a corollary of Theorem \ref{thm-main-intro}, at least if $G$ is a connected and simply connected Lie group. 

Since actions $\beta:L\car A$ of compact groups $L$ on type I $C^*$-algebras $A$ are always smooth and the crossed products $A\rtimes_\beta L$ are then also type I, another corollary of Theorem \ref{thm-main-intro} is the following

\begin{corollary}\label{cor-compact-group}
Suppose that $\beta:L\car A$ is an action of a {\em compact} abelian group $L$ on the type I $C^*$-algebra 
$A$ and let $\varphi:G\car L$ be a continuous homomorphism from the abelian locally compact group $G$ into $L$.
Let $\alpha=\beta\circ \varphi:G\car A$ be the pull-back of $\beta$ via $\varphi$. Then the conclusion of
Theorem \ref{thm-main-intro} holds for $A\rtimes_\alpha G$.
\end{corollary}

Note that the action $\alpha:G\car A$ in the above corollary can be quite far away from being smooth itself. Indeed, interesting examples for the corollary come from inclusions $\varphi:G\into L$ of infinite countable groups, in which case the action $\alpha=\beta\circ \varphi$ will never be smooth!

Our proof of Theorem \ref{thm-main-intro} is of course inspired by the ideas presented by Poguntke in his paper \cite{Pogi}.
But we tried (and hopefully succeeded) to avoid some of the tedious and technical
arguments given in the original proof by using some general concepts 
in the realm of the Mackey-Rieffel-Green machine, 
which partly  have been developed 
after the publication of Puguntke's paper \cite{Pogi}. 
 So we hope that our exposition will contribute to a better understanding of this fundamental result.

These notes are outlined as follows: we start in Section \ref{sec-Mackey-Rieffel-Green} with a brief account on twisted crossed products and some useful applications of Green's imprimitivity theorem.  
In Section \ref{sec-subquotients} we study the decomposition of crossed products by subquotients as mentioned above. In Section \ref{sec-actionsK}, we explain the construction of twisted group algebras $C^*(G,\om)$ 
for a Borel cocycle $\om\in Z^2(G,\T)$ and explain the relation between twisted group algebras and crossed products by 
actions $\alpha:G\car \K(\H)$ on the algebra of compact operators $\K(\H)$. 
In  Section \ref{sec-connected} we present the proofs of Theorem \ref{thm-main-intro} and Corollary \ref{cor-compact-group} as stated above. 

In Section \ref{sec-twisted-group-algebras} we  study of the  structure of twisted group algebras 
$C^*(G,\om)$ if $G$ is abelian. 
Based on the work of Baggett and Kleppner, we give a proof of the (well-known) result
that  $C^*(G,\om)$ is a continuous field of $C^*$-algebras over the Ponrjagin dual 
$\widehat{S_\om}$ of the 
symmetry group $S_\om$, with fibres isomorphic to 
a simple twisted group algebra $C^*(G/S_\om, \tilde\om)$ for a totally skew cocycle $\tilde\om\in Z^2(G/S_\om,\T)$ that inflates to $\om$ on $G$. We then give several reduction steps on the structure of the group $G$ for twisted group algebras by totally skew cocycles up to stable isomorphism. In particular we show  that simple twisted group algebras 
of compactly generated abelian groups are elementary or Morita equivalent to a non-commutative torus.
This result is equivalent to Poguntke's \cite{Pogi}*{Theorem 1} and we do not claim any originality for the approach taken here.

In  Section \ref{sec-Poguntke} we obtain  Corollary \ref{cor-main-intro}  as a combination of 
Theorem \ref{thm-main} and the results of Section \ref{sec-twisted-group-algebras}. We  then use this corollary to 
show that the simple subquotients of group algebras of connected groups are either compact operators or Morita equivalent to simple non-commutative tori of dimension $\geq 2$. Together with the results of Section \ref{sec-twisted-group-algebras} this gives the proof of Theorem \ref{thm-Pogi}. In the  final Section \ref{sec-almost connected} we briefly study the case of almost cennected groups $G$, i.e., groups with co-compact connected component $G_0$. As we shall see, the 
primitive ideals are still all locally closed and each non-elementary simple subquotient of $C^*(G)$ is Morita equivalent to a 
quotient of some crossed product $(C^*(\Z^m,\om)\otimes \K(\ell^2(\N)))\rtimes_\gamma F$ with  $C^*(\Z^m,\om)$ a simple non-commutative torus and  $F$ finite. The general structure of such crossed products is still a mystery, but  at least for actions $\gamma:F\car C^*(\Z^m,\om)$ some partial results have been obtained in the literature (e.g., see \cites{ELPW, Chakraborty, Jeong, He}).

\subsection*{Aknowlegements}  The original motivation for this work came from a lecture by Caleb Eckhard on the classifiability of 
simple quotients of $C^*$-group algebras of virtually nilpotent groups (see \cite{Eckhardt}). In the light of Poguntke's results, it seemed to be natural to assume that some of the methods presented by Caleb could possibly also lead to a better understanding of simple subquotients of $C^*$-group algebras of almost connected groups. As a starting point for such a research, it seemed necessary to get a deeper understanding of Poguntke's results and their proofs. We are very grateful to Caleb Eckhard and Ilan Hirshberg for the many motivating and fruitful discussions we had on this topic.

\section{Twisted crossed products and Morita equivalence}\label{sec-Mackey-Rieffel-Green}
In this section we want to recall some basics on the Mackey-Rieffel-Green machine for crossed products by actions 
$\alpha:G\car A$ of a locally compact group $G$ on a $C^*$-algebra $A$. For more complete accounts we refer 
to the original paper  by Phil Green \cite{Green},  the book on crossed-products by Dana Williams \cite{Dana:book},
or the survey given in \cite{CELY}*{Chapter 2}. The latter reference probably fits best to the methods used here.

We start with a brief recollection of the construction of twisted crossed products in the sense of Green \cite{Green}.
For this let $G$ be a locally compact group and let $N\subseteq G$ be a closed  normal subgroup.
A (Green-)\emph{twisted} action $(\alpha, \tau)$ of  the pair $(G,N)$ on a C*-algebra $A$ consists 
of a (continuous) action $\alpha:G\car A$ by $*$-automorphisms together with a strictly continuous homomorphism 
$\tau:N\to U\M(A)$ such that
\begin{equation}\label{eq-Green-twist}
\alpha_g(\tau_n)=\tau_{gng^{-1}}\quad\text{and}\quad \alpha_n(a)=\tau_n a \tau_{n^{-1}}
\end{equation}
for all $g\in G, n\in N$, and $a\in A$.

If $(\alpha,\tau)$ is a twisted action as above and if $(\pi, U)$ is a non-degenerate covariant representation\footnote{So, for all $g\in G, a\in A$, it satisfies the equation $\pi(\alpha_g)=U_g\pi(a)U_g^*$.}  of the underlying 
system $(A,G,\alpha)$ on a Hilbert space (or Hilbert module) $\H$, or into the multiplier algebra $\M(C)$ of some  $C^*$-algebra 
$C$, say, then we say that $(\pi, U)$ {\em preserves the twist $\tau$} if for all $n\in N$:
\begin{equation}\label{eq-Green-twist1}
\pi(\tau_n)=U_n.
\end{equation}
We then say that $(\pi,U)$ is a {\em covariant representation of  the twisted system} \linebreak
$(A,G, N,\alpha,\tau)$.
Note that since we required $\pi$ to be non-degenerate, it extends uniquely to $U\M(A)$, and hence the
equation makes sense. 
\medskip

The (Green-) {\em twisted crossed product} $A\rtimes_{(\alpha,\tau)}G/N$ can then be defined as the quotient 
of the ordinary crossed product $A\rtimes_\alpha G$ by the {\em twisting ideal} 
$$I_{\tau}=\cap\{\ker(\pi\rtimes U): (\pi, U)\;\text{preserves $\tau$}\},$$
where $\pi\rtimes U$ denotes the usual integrated form of $(\pi,U)$. 
\medskip 

Another possible construction of the twisted crossed product is by considering the 
space
\begin{equation}\label{eq-Green-twist2}
C_c(G, A,\tau):=\left\{f:G\to A: \begin{matrix} \text{$f$ cont.,  $f(ng)=f(g)\tau_{n^{-1}}$ $\forall n\in N, g\in G$,}\\
\text{and $(gN\mapsto \|f(g)\|)\in C_c(G/N)$}\end{matrix}\right\}
\end{equation}
equipped with convolution and involution given by
\begin{equation}\label{eq-Green-twist3}
f_1*f_2(g)=\int_{G/N} f_1(h)\alpha_h(f_2(h^{-1}g))\, d (hN) \quad\text{and}\quad f^*(g)=\Delta_{G/N}(g^{-1})\alpha_g(f(g^{-1})^*),
\end{equation}
for $f, f_1, f_2\in C_c(G, A,\tau)$ and $g\in G$. 
If $(\pi, U)$ is a covariant representation of $(A,G,N, \alpha, \tau)$
we consider the integrated form
$$\pi\rtimes U(f):=\int_{G/N} \pi(f(g))U_g\, d gN,\quad f\in C_c(G,A,\tau)$$
of $(\pi,U)$  and then  $A\rtimes_{(\alpha,\tau)}G/N$ can also be realized as the completion of $C_c(G,A,\tau)$ 
with respect to 
$$\|f\|:=\sup\{\|\pi\rtimes U(f)\|: \text{$(\pi,U)$ a covariant rep of $(A, G, N, \alpha, \tau)$}\}.$$
\medskip
 
A twisted action of the pair $(G,N)$ should be regarded as a generalization of an ordinary action of the quotient group $G/N$, hence our notation 
for the crossed product.
To get a hint for this, let us first look at the case $N=G$. If $(\alpha,\tau):(N,N)\car A$ is a twisted action, we get
$A\rtimes_{(\alpha,\tau)}N/N\cong A$ via
$$A\to A\rtimes_{(\alpha,\tau)}N/N; a\mapsto (n\mapsto a\tau_{n^{-1}})\in C_c(N, A,\tau).$$
So twisted actions by the pair $(N,N)$ behave like actions of the trivial group $\{e\}$.
Conversely, the following example shows that every action of $G/N$ can be regarded as a twisted action of 
$(G,N)$ in a canonical way:

\begin{example}\label{ex-inflate}[(cf \cite{Ech:Morita}*{p.~176})]
Suppose that $\alpha:G/N\car A$ is an ordinary action of the quotient group $G/N$ on $A$.  Let $\tilde\alpha:G\car A$ denote the 
{\em inflation} of $\alpha$ defined by $\tilde\alpha:=\alpha\circ q$, where $q:G\to G/N$ denotes the quotient map.
Moreover, let $1_N:N\to U\M(A); n\mapsto 1_{\M(A)}$ denote the trivial homomorphism. 
Then $(\tilde\alpha, 1_N)$ is a twisted action of $(G,N)$ on $A$ such that 
$A\rtimes_{(\tilde\alpha,1_N)}G/N\cong A\rtimes_\alpha G/N$. 
Indeed, we directly check that $C_c(G,A,1_N)\cong C_c(G/N, A)$, since the condition on $f\in C_c(G,A,\tau)$  in (\ref{eq-Green-twist3})
for $\tau=1_N$ implies that $f$ is constant on $N$-cosets. Moreover, a covariant representation 
$(\pi, U)$ of $(A, G, \tilde\alpha)$ preserves $1_N$ as in (\ref{eq-Green-twist1}) if and only if 
$U$ factors through a representation $\dot{U}$ of $G/N$, and hence corresponds to the 
covariant representation $(\pi, \dot{U})$ of $(A,G/N, \alpha)$ (and vice versa). 
The claim then follows easily from the constructions.
\end{example}

But the 
motivating example for the introduction of twisted actions and crossed products as above is the following construction due to Green
(\cite{Green}):

\begin{example}\label{ex-decomp}[(cf \cite{Green}*{Proposition 1})]
Suppose that $\beta:G\car B$ is any ordinary action of $G$ on a $C^*$-algebra $B$ and let $N$ be a closed normal subgroup 
of $G$. We want to write the crossed product 
$B\rtimes_\beta G$ as an iterated crossed product $(B\rtimes_\beta N)\rtimes_{\tilde\beta} G/N$.
This makes perfect sense if $G=N\rtimes H$ is a semi-direct product (hence $H\cong G/N$).  But for non-split extensions this is not possible 
within the  realm of ordinary actions and crossed products, since in general there does not exist a well defined action of the quotient group 
$G/N$ on $B\rtimes_{\beta}N$. But we do have a well defined twisted action $(\tilde\beta, \sigma)$ of the pair $(G,N)$
on $B\rtimes_\beta N$ given on the level of the dense subalgebra $C_c(N,B)$ as follows:
\begin{equation}\label{eq-Green-twist4}
\big(\tilde\beta_g(f)\big)(n)=\frac{\Delta_G(g)}{\Delta_{G/N}(g)}\beta_g(f(g^{-1}ng))\quad\text{and}\quad \big(\sigma_m f\big)(n)=\beta_m(f(m^{-1}n)).
\end{equation}
for $f\in C_c(N,B)$ and $m,n\in N$. 
The mapping 
$$\Phi:C_c(G,B)\mapsto C_c(G, C_c(N,B), \sigma),\quad \Phi(f)(g)(n)=\frac{\Delta_G(g)}{\Delta_{G/N}(g)}f(ng)$$
then extends to an isomorphism
$$ B\rtimes_\beta G\cong (B\rtimes_\beta N)\rtimes_{(\tilde\beta,\sigma)}G/N.$$
Similarly, if we start with a twisted action $(\beta,\tau): (G,M)\car B$ and $N$ is a normal subgroup of $G$ such that 
$M\subset N$, then we obtain a twisted action of $(\tilde\beta, \sigma): G/N\car B\rtimes_{(\beta,\tau)}N/M$ such that 
$$B\rtimes_{(\beta,\tau)}G/M\cong \big(B\rtimes_{(\beta,\tau)}N/M\big)\rtimes_{(\tilde\beta, \sigma)}G/N,$$
where in the above description of the decomposition action we may replace $C_c(N,B)$ by $C_c(N,B,\tau)$.
\end{example}

We proceed  with a discussion on  induced $C^*$-algebras and Green's imprimitivity theorem.
Suppose that $H$ is a closed subgroup of the locally compact group $G$ and that
$\beta:H\car B$ is an action of $H$ on the C*-algebra $B$ by $*$-automorphisms. Then the {\em induced $G$-C*-algebra}
$\Ind_H^G(B,\beta)$ (or just $\Ind_H^GB$ if $\beta$ is understood) is defined as the set of $B$-valued functions
\begin{equation}\label{eq-induced}
\Ind_H^G(B,\beta)=\left\{f:G\to B: \begin{matrix} \text{$f$ continuous, $f(gh)=\beta_{h^{-1}}(f(g))$  $\forall g\in G, h\in H$, }\\
\text{and $(gH\mapsto \|f(g)\|)$ vanishes at $\infty$ on $G/H$.}\end{matrix}\right\}
\end{equation}
equipped with pointwise multiplication and involution,  and the supremum norm. The {\em induced action} 
of $G$ on $\Ind_H^GB$ is given by
\begin{equation}\label{eq-induced-action}
(\Ind\beta_g(f)(l)=f(g^{-1}l)\quad \forall g,l\in G.
\end{equation}
More generally, if $N\subseteq G$ is a closed normal subgroup contained in $H$ and $\sigma:N\to U\M(B)$ is a twist for $\beta$,
then we obtain a twist $\Ind\sigma:N\to U\M(\Ind_H^GB)$ for $\Ind\beta$ by putting $(\Ind\sigma(n)f)(g)=f(n^{-1}g)$ ($=\sigma_{g^{-1}ng}f(g)$).\footnote{Observe that 
$\Ind\sigma_n:g\mapsto \sigma_{g^{-1}ng}$ is a function $G\to U\M(B)$ that satisfies the condition $\Ind\sigma_n(gh)=\beta_{h^{-1}}(\Ind\sigma(g))$ for all $g\in G, h\in H$,
and therefore defines an element in $U\M(\Ind_H^GB))$ via point-wise multiplication.}
This provides 
us with the 
{\em induced twisted action} $(\Ind\beta,\Ind\sigma):(G,N)\car \Ind_H^GB$.

\begin{remark}\label{rem-trivial action}
 If $\beta:H\car B$ extends to an action $\tilde\beta:G\car B$, then the induced system 
$(\Ind_H^GB, G, \Ind\beta)$ is isomorphic to the system $(C_0(G/H)\otimes B, G, \lt\otimes\tilde\beta)$,
where $\lt:G\car C_0(G/H)$ is the action by left translation. The isomorphism is given by
$$\Phi: C_0(G/H, B)\to \Ind_H^GB; \quad \Phi(f)(g)=\tilde\beta_{g^{-1}}(f(gH)).$$
If $\sigma:N\to U\M(B)$ is a twist for $\beta$, the induced twist realized on $C_0(G/H,B)=C_0(G/H)\otimes B$ is 
simply given by $1\otimes \sigma$. 
\end{remark}

Observe that the primitive ideal space $\Prim(\Ind_H^GB)$ can be identified with the set
$$G\times_H\Prim(B):= (G\times \Prim(B))/H$$
with action $H\car G\times \Prim(B)$ given by $h\cdot(g, P)=(gh^{-1}, \beta_h(B))$. 
The primitive ideal $P_g$ of $\Ind_H^GB$ corresponding to an equivalence class $[g, P]$ is given by
$$P_g:=\{f\in \Ind_H^GB: f(g)\in P\}.$$
Notice that  $\varphi: \Prim(\Ind_H^GB)\to G/H; \varphi(P_g)=gH$ is a continuous $G$-equivariant map.
This observation lead to the following  general characterization of induced algebras:

\begin{theorem}[{\cite{Ech:induced}*{Theorem}}]\label{thm-induced}
Let $\alpha:G\car A$ be an action of the locally compact group $G$ on the C*-algebra $A$ and let $H$ be a closed subgroup of $G$.
Then the system $(A,G,\alpha)$ is equivariantly isomorphic to an induced system $(\Ind_H^GB, G,\Ind\beta)$ for some action 
$\beta:H\car B$ if and only if there exists a $G$-equivariant continuous mapping $\varphi:\Prim(A)\to G/H$. 

If, moreover,  $\tau:N\to U\M(A)$ is a twist for $\alpha$ such that $N\subseteq H$,
then there exists a twist $\sigma:N\to U\M(B)$ such that 
$(A,G,N, \alpha,\tau)$ is isomorphic to the induced twisted system $(\Ind_H^GB, G,\Ind\beta,\Ind\sigma)$.
\end{theorem}

Indeed,  if $\varphi:\Prim(A)\to G/H$ is as in the theorem, we can choose $B:=A/I$ with $I:=\cap\{P\in \Prim(A): \varphi(P)=eH\}$.
This is an $H$-invariant quotient of $A$, and therefore the restriction of $\alpha$ to $H$ induces a well-defined action 
$\beta:H\car B$. The isomorphism  $\Phi:A\to \Ind_H^GB$ is then given by sending an element $a\in A$ to the function 
$f_a\in \Ind_H^GB$ given by
$$f_a(g)=(\alpha_{g^{-1}}(a)+I)\in B.$$
It is shown in \cite{Ech:induced}*{Theorem} that this is an isomorphism for the untwisted systems. 
Now, if we define $\sigma_n$ as the image of $\tau_n$ under the extension to $\M(A)$  of the quotient map
$A\to A/I=B$, then for all $a\in A$ we get
$$\Phi(f_{\tau_n a})(g)=(\alpha_{g^{-1}}(\tau_n a)+I)=\big(\big(\tau_{g^{-1}ng}\alpha_{g^{-1}}(a)\big)+I\big)=\big(\Ind\sigma_n\Phi(f)\big)(g),$$
so we see that $\Phi$ also preserves the twists.

With these notations, Green's imprimitivity theorem (see \cite{Green}*{Theorem 17 and Theorem 18}) can be formulated as follows:

\begin{theorem}[Green]\label{thm-Green}
Suppose $H$ is a closed subgroup of $G$  and let $\beta:H\car B$ be an action. Then 
$B\rtimes_\beta H$ is Morita equivalent to $\Ind_H^GB\rtimes_{\Ind\beta}G$.
A similar result holds for an induced twisted action $(\Ind\beta,\Ind\sigma):(G,N)\car \Ind_H^GB$.
\end{theorem}

 Although, the above theorem is covered by (and indeed equivalent to) \cite{Green}*{Theorem 18}, 
it is important for us to recall a more direct  construction of an  $\Ind_H^GB\rtimes_{\Ind\beta}G-B\rtimes_\beta H$ equivalence bimodule 
in the situation of the above theorem as given in \cite{Dana:book}*{Corollary 4.17}. 
There it is realized  as the completion $\mathcal X_H^G(B)$ of the $C_c(G,\Ind_H^GB)-C_c(H,B)$ bimodule $\mathcal X_0:=C_c(G,B)$
with inner products and actions given by the formulas
\begin{equation}\label{eq-imp}
\begin{split}
f\cdot \xi(g)&= \int_G f(s, g) \xi(s^{-1}g) \Delta_G(s)^{1/2}\, d s\\
\xi\cdot \varphi(g)&= \int_H \beta_{h^{-1}}\big(\xi(gh^{-1})\varphi(h)\big)\Delta_H(h)^{-1/2} \, dh\\
_{C_0(G,\Ind_H^GB)}\lk \xi, \eta\rk(g,s)&= \Delta_G(g)^{-1/2}\int_H \beta_h\big(\xi(sh)\eta(g^{-1}sh)^*\big) \, dh\\
\lk \xi,\eta\rk_{C_c(H,B)}(h)&=\Delta_H(h)^{-1/2}\int_G \xi(g^{-1})^*\beta_h(\eta(g^{-1}h))\, dg
\end{split}
\end{equation}
for $f\in C_c(G,\Ind_H^GB), \varphi\in C_c(H,B)$, and $\xi,\eta\in \mathcal X_0$.
This provides the desired Morita equivalence in the untwisted situation. In the twisted case, it suffices to check that
a covariant representation $\rho\times V$ of $B\rtimes_\beta H$ preserves the twist $\sigma$ if and only if the induced 
representation $\Ind^{\X_H^G(B)}(\rho\times V)$ of $\Ind_H^GB\rtimes_{\Ind\beta}G$ preserves the twist $\Ind\sigma$.
We leave the verification as an exercise for the interested reader.

\begin{remark}\label{rem-induced-rep}
Note that Green's original imprimitivity theorem (\cite{Green}*{Theorem 17}) was stated for the special case where $\beta:H\car B$ is the 
restriction of  an action $\beta:G\car B$, i.e., where $(\Ind_H^GB,\Ind\beta)$ is equivariantly isomorphic to $(C_0(G/H)\otimes B, \lt\otimes\beta)$ (see Remark \ref{rem-trivial action} above). The resulting $C_0(G/H,B)\rtimes_{\lt\otimes\beta}G-B\rtimes_\beta H$ equivalence bimodule $\X_H^G(B)$   provides a convenient way
to define induction of covariant representations from $H$ to $G$:
given a non-degenerate representation $\rho\rtimes V$ of $B\rtimes_\beta H$, we can first induce it to $(B\otimes C_0(G/H))\rtimes_{\beta\otimes \lt}G$ via the equivalence 
bimodule $\X_H^G(B)$ (as introduced by Rieffel in \cite{Rief:InducedC}), and then compose it with the $*$-homomorphism
$$\iota_B\rtimes G: B\rtimes_\beta G\to \M((C_0(G/H,B)\rtimes_{\lt\otimes \beta} G),$$
where $\iota_B:B\to \M(C_0(G/H, B))$ denotes the canonical inclusion by constant functions.  In short, the induced representation $\Ind_H^G(\rho\rtimes V)$ is defined as
$$\Ind_H^G(\rho\rtimes V):=\big(\Ind^{\X_H^G(B)}(\rho\times V)\big)\circ (\iota_B\rtimes G).$$
Exactly the same works in the twisted case.
It is shown in \cite{Ech:T1} that this definition coincides with the more classical definition due to Takesaki using  constructions 
of  Mackey and  Blattner.
It is trivial to deduce from Green's description, that induction of representations is continuous in the sense that it preserves weak containment of representations in the sense of Fell.
\end{remark}

Allmost all  results one can prove for ordinary actions can be shown equally well 
for twisted actions and their crossed products. However, many formulations and proofs become a bit more clumsy by carrying around the twists all the time.
A nice way out is given by the Packer-Raeburn stabilization theorem \cite{PackRae}*{Theorem 3.4} which, in case $G$ is second countable and $A$ is separable, provides an ordinary 
action $\beta:G/N\car A\otimes \K$, $\K=\K(\ell^2(\N))$, such that 
$$(A\otimes \K)\rtimes_\beta G\cong \big(A\rtimes_{(\alpha,\tau)}G/N\big)\otimes \K.$$
A slightly modified version, shown in \cite{Ech:Morita}, works without any separability assumptions and enjoys some better 
functoriallity conditions. We need

\begin{definition}[\cite{Ech:Morita}]\label{def-Mor}
We say that two twisted actions $(\alpha,\tau):G/N\car A$ and $(\beta,\sigma):G/N\car B$ are {\em Morita equivalent}
if there exists an $A-B$ equivalence bimodule $\mathcal{X}$ together with an action $\gamma:G\car \mathcal{X}$ such that
\begin{enumerate} 
\item $(\mathcal{X},\gamma)$ induces an $(A,\alpha)-(B,\beta)$ Morita equivalence in the sense that
\begin{align*} _A\lk \gamma_g(x), \gamma_g(y)\rk&=\alpha_g(_A\lk x,y\rk),\quad  \lk \gamma_g(x), \gamma_g(y)\rk_B=\beta_g(\lk x,y\rk_B)\\
\text{and}\;& \gamma_g(ax)=\alpha_g(a)\gamma_g(x), \quad \gamma_g(xb)=\gamma_g(x)\beta_g(b)
\end{align*}
for all $x,y\in \mathcal{X}, a\in A, b\in B$, and $g\in G$.
\item $\tau_n x=\gamma_n(x)\sigma_n$ for all $n\in N, x\in \mathcal{X}$.
\end{enumerate}
If $(\alpha,\tau):(G,N)\car A$ is Morita equivalent to a twisted action $(\tilde\beta, 1_N):(G,N)\car B$ inflated from an 
ordinary action $\beta:G/N\car B$, we shall say that $(\alpha,\tau)$ is Morita equivalent to the (ordinary) action $\beta$.
\end{definition}

\begin{remark}\label{rem-twistedMor}
{\bf (1)} If $(\X,\gamma)$ is a Morita equivalence for the twisted actions $(\alpha,\tau):(G,N)\car A$ and $(\beta,\sigma):(G,N)\car B$, then 
we obtain a twisted action 
$$(\eta,\mu):=\left(\left(\begin{matrix}\alpha&\gamma\\ \gamma^*&\beta\end{matrix}\right),\left(\begin{matrix} \tau&0\\ 0&\sigma\end{matrix}\right)\right)
:(G,N)\car L(\mathcal{\X})$$
on the linking algebra $L(\mathcal{\X}):=\left(\begin{smallmatrix} A&\mathcal X\\ \mathcal X^*& B\end{smallmatrix}\right)$. 
Conversely, every action $\gamma:G\car \X$ that fits into the above matrix to give a twisted action of $(G,N)$ on $\L(\mathcal X)$ induces an
$(\alpha,\tau)-(\beta,\sigma)$ Morita equivalence. The  crossed product $L(\mathcal{X})\rtimes_{(\eta,\mu)}G/N$ then naturally decomposes
into a linking algebra
$$L(\mathcal{X})\rtimes_{(\eta,\mu)}G/N=:\left(\begin{matrix} A\rtimes_{(\alpha,\tau)}G/N& \X\rtimes_\gamma G/N\\ \X^*\rtimes_{\gamma^*}G/N& B\rtimes_{(\beta,\sigma)}G/N
\end{matrix}\right)$$
for an $A\rtimes_{(\alpha,\tau)}G/N-B\rtimes_{(\beta,\sigma)}G/N$ equivalence bimodule $\X\rtimes_{\gamma}G/N$. 
Thus, Morita equivalent twisted actions have Morita equivalent twisted crossed products. 

{\bf (2)} Recall from \cite{Rief:InducedC} that, by induction via $\mathcal X$,  any $A-B$ Morita equivalence $\X$ induces a homeomorphism $\widehat{A}\cong\widehat{B}$ of the dual spaces, 
and a homeomorphism  $\Prim(A)\cong \Prim(B)$ between the primitive ideal spaces, equipped with the Jacobson topologies.
If $(\X,\gamma)$ is an equivariant Morita equivalence between the twisted actions $(\alpha,\tau):(G,N)\car A$ and $(\beta,\sigma):(G,N)\car B$, then these 
homeomorphisms are $G$-equivariant.

{\bf (3)} 
If $G/N$ is abelian,  the dual action $\widehat{(\alpha,\tau)}:\widehat{G/N}\car A\rtimes_{(\alpha,\tau)}G/N$ is defined by 
\begin{equation}\label{eq-dualaction}
\widehat{(\alpha,\tau)}_\chi(f)=\bar\chi\cdot f
\end{equation}
for $\chi\in \widehat{G/N}$ and $f\in C_c(G, A,\tau)$. Now, if $(\alpha,\tau)$ is Morita equivalent to $(\beta, \sigma)$ via $(\mathcal{X}, \gamma)$,
the  dual action of $\widehat{G/N}$ on 
$$L(\mathcal{X})\rtimes_{(\eta,\mu)}G/N=:\left(\begin{matrix} A\rtimes_{(\alpha,\tau)}G/N& \X\rtimes_\gamma G/N\\ \X^*\rtimes_{\gamma^*}G/N& B\rtimes_{(\beta,\sigma)}G/N
\end{matrix}\right)$$
restricts to a  dual action $\widehat\gamma:\widehat{G/N}\car \mathcal{X}\rtimes_\gamma G/N$ implementing  a Morita equivalence between the dual actions 
$\widehat{(\alpha,\tau)}$ and $\widehat{(\beta,\sigma)}$. 
\end{remark}

The above observations  shows that it is  very helpful to know 
that {\em every} twisted action $(\alpha,\tau):(G,N)\car A$ is Morita equivalent to some 
ordinary action $\beta:G/N\car B$ as above. This will then allow us to pass to
ordinary actions and crossed products in all situations where the desired results only depend on  
Morita equivalence classes (like most results we want to discuss in this paper).

\begin{theorem}[\cite{Ech:Morita}]\label{thm-Morita}
Suppose that $(\alpha,\tau):G/N\car A$ is a twisted action of $(G,N)$ on $A$. Consider the 
twisted crossed product 
$$B:=(C_0(G/N)\otimes A)\rtimes_{(\lt\otimes \alpha, 1\otimes \tau)}G/N,$$
where $\lt:G\car C_0(G/N)$ denotes the action by left translation.
Then the right translation action  $\rt:G/N\car C_0(G/N)$ commutes with $(\lt\otimes \alpha\otimes \lt, 1\otimes \tau)$ and therefore induces an action 
$\beta:G/N\car B$. Then $(\alpha, \tau)$ is Morita equivalent to $\beta$.
\end{theorem}

Recall that the algebra $B$ above is Green's imprimitivity algebra for induction from $A\rtimes_{(\alpha,\tau)}N/N\cong A$ to 
$A\rtimes_{(\alpha,\tau)}G/N$. Then Green's impriimitivity theorem provides a
$B-A$ equivalence bimodule $\X_N^G(A)$ and the Morita equivalence of the above theorem is given by 
a suitable action $\gamma:G\car \X_N^G(A)$. We refer to \cite{Ech:Morita}*{Theorem 1}  for the details.

\begin{example}\label{ex-abelian}
If $G/N$ is abelian, we can consider the dual action $\widehat{(\alpha,\tau)}:\widehat{G/N}\car A\rtimes_{(\alpha,\tau)}G/N$ as in Remark \ref{rem-twistedMor}
above. Dualising again,  we may consider the double dual action 
$$\widehat{\widehat{(\alpha,\tau)}}:G/N\car A\rtimes_{(\alpha,\tau)}G/N\rtimes_{\widehat{(\alpha,\tau)}}\widehat{G/N}.$$
It is shown in \cite{Ech:Duality} that this action is indeed isomorphic to the action $\beta:G/N\car B$  of Theorem \ref{thm-Morita},
and hence we see that $(\alpha,\tau)$  is Morita equivalent to its double dual action. Note that for ordinary actions 
$\alpha:G\car A$ of abelian groups $G$, the Morita equivalence $\alpha\sim_M\widehat{\widehat{\alpha}}$ is  a direct 
consequence of the Takesaki-Takai duality theorem, which identifies  $A\rtimes_\alpha G\rtimes_{\widehat\alpha}\widehat{G}$ 
with $A\otimes \K(L^2(G))$ and $\widehat{\widehat{\alpha}}$ with $\alpha\otimes \Ad \rho$, where $\rho:G\to \U(L^2(G))$ 
denotes the right regular representation.
\end{example}

\begin{example}\label{ex-induced-Mor}
Let $H$ be a closed subgroup of the abelian locally compact group $G$ and let $\alpha:H\car A$ be an action of $H$ 
on the $C^*$-algebra $A$. Then the inflation $\Inf\widehat\alpha:\widehat{G}\car A\rtimes_\alpha H$ of the dual action $\widehat\alpha:\widehat{H}\car A\rtimes_\alpha H$ given by
$\Inf\widehat\alpha_\chi=\widehat\alpha_{\chi|_H}$
is Morita equivalent to the dual action $\widehat{\Ind\alpha}:\widehat{G}\car \Ind_H^GA\rtimes_{\Ind\alpha}G$.
Indeed, it follows easily from the 
formulas given in (\ref{eq-imp}) that the action $\gamma: \widehat{G}\car \mathcal X_H^G(A)$ of $\widehat{G}$ 
on Green's  $\Ind_H^GA\rtimes_{\Ind\alpha}G-A\rtimes_\alpha H$ equivalence bimodule $\mathcal X_H^G(A)$ given on the 
dense submodule $\mathcal X_0$ by 
\begin{equation}\label{eq-Ghat-equivariance}
\gamma_\chi(\xi)(g)=\overline{\chi(g)}\xi(g),\quad \xi\in \mathcal X_0
\end{equation}
implements the desired equivalence for the actions. 
\end{example}

\section{Decomposition by subquotients}\label{sec-subquotients}

Recall from the introduction that a subset $Y$ of the primitive ideal space $\Prim(A)$ of a $C^*$-algebra $A$ 
is locally closed 
in $\Prim(A)$  if and only if there 
exist closed ideals $J\subseteq I\subseteq A$ of $A$, such that $Y\cong \Prim(I/J)$.
The ideals $I$ and $J$ are uniquely determined by $Y$ as 
$$J=\ker (Y)\quad \text{and} \quad I=\ker({\overline{Y}\setminus Y}),
$$
and we  then write $A_Y:=I/J$.
A similar correspondence holds between locally closed subsets of  the space $\widehat{A}$ of equivalence classes of irreducible $*$-representations of $A$ and subquotients of $A$.
For later use we need the following proposition, which is certainly  well-known to experts.
We leave the straightforward proof as an exercise for the reader.

\begin{proposition}\label{prop-loc-closed-points}
Let $X$ be a topological space and let $Z\subseteq Y\subseteq X$. 
If $Z$ is locally closed in $X$, it is locally closed in $Y$.
Conversely, if $Y$ is locally closed in $X$ and $Z$ is
locally closed in $Y$, then $Z$ is locally closed in $X$.

Moreover, if 
$Z\subseteq Y$ are locally closed subsets of  $X=\Prim(A)$ for some $C^*$-algebra $A$, and 
$I_Z, J_Z$ and $I_Y, J_Y$ are closed ideals in $A$ such that 
$$Z\cong \Prim(I_Z/J_Z)\quad\text{and}\quad Y\cong \Prim(I_Y/J_Y),$$
and  $I_Z^B, J_Z^B$ are the  ideals in $B:=I_Y/J_Y$ corresponding to 
$Z\subseteq Y=\Prim(B)$, then $I_Z^B/J_Z^B\cong I_Z/J_Z$.
\end{proposition}

In what follows next, we want to study decompositions by subquotients of crossed products $A\rtimes_\alpha G$ 
coming from $G$-invariant decompositions of the algebra $A$. We restrict our attention to the case of untwisted crossed 
products, noting that the twisted case can easily be  achieved from  the untwisted case via Theorem \ref{thm-Morita}.

If $\alpha:G\car A$ is an action of $G$ on the $C^*$-algebra $A$ and if $Y\subseteq \Prim(A)$ is locally closed 
and $G$-equivariant, then the ideals $J=\ker Y$ and $I=\ker(\overline{Y}\setminus Y)$ are $G$-invariant as well. Hence 
the  action $\alpha$ induces   actions $\alpha^J:G\car J$ and $\alpha^I:G\car I$  by restriction, 
and then  an action $\alpha^Y:G\car A_Y$ on the quotient $A_Y=I/J$ (we shall often simply write $\alpha:G\car A_Y$ if confusion seems unlikely). Since taking full crossed products sends short exact sequences of 
$G$-algebras to short exact sequences of $C^*$-algebras, 
we see that
$$A_Y\rtimes_{\alpha^Y}G\cong \big(I\rtimes_{\alpha^I}G\big)/ \big(J\rtimes_{\alpha^J}G\big)$$
is a subquotient of $A\rtimes_{\alpha}G$, and therefore we can conclude that $\Prim\big(A_Y\rtimes_{\alpha^Y}G\big)$
is a locally closed subset of $\Prim(A\rtimes_{\alpha}G)$. 
The proof of the following proposition is modelled after \cite{Green}*{Corollary 19}. It easily follows from Theorem \ref{thm-Morita} 
that a  similar result holds true for twisted crossed products as well.

\begin{proposition}\label{prop-orbit-method}
 Suppose that $\alpha:G\car A$ is an action of the locally compact group $G$ on the $C^*$-algebra $A$ and that  $\phi:\Prim(A)\to X$ 
 is a $G$-invariant continuous open map onto the almost Hausdorff topological space $X$.
 Then $\phi^{-1}(x)$ is locally closed in $\Prim(A)$ for all $x\in X$, and  if $A_x:=A_{\phi^{-1}(x)}$ is the subquotient of $A$
 corresponding to $\phi^{-1}(x)$, then the  action of $G$ on $A$ implements actions $\alpha_x:G\car A_x$ such that 
 $(A_x)_{x\in X}$ is a decomposition of $A$ and $(A_x\rtimes_{\alpha_x}G)_{x\in X}$ is a decomposition of $A\rtimes_\alpha G$ by 
 subquotients. 
\end{proposition}

\begin{proof}  Since $\Prim(A)$ is the disjoint union of the locally closed sets $\phi^{-1}(x)$, $x\in X$, the assertion on $(A_x)_{x\in X}$ 
is clear. To see that $(A_x\rtimes_{\alpha_x}G)_{x\in X}$ is a decomposition of $A\rtimes_\alpha G$, we need to show that 
for every primitive ideal $Q\in \Prim(A\rtimes_\alpha G)$ there exists a unique $x\in X$ such that $Q$ belongs to
$\Prim(A_x\rtimes_{\alpha_x}G)$. 

So let $Q\in \Prim(A\rtimes_\alpha G)$ be fixed and let $\pi\rtimes U$ be an irreducible representation of $A\rtimes_\alpha G$ 
such that $Q=\ker(\pi\rtimes U)$. For each $x\in X$ let $J_x\subseteq I_x\subseteq A$ be the closed ideals such that 
$A_x=I_x/J_x$.  Since the kernel $\ker \pi\subseteq A$ only depends on $Q$ and not  on the choice of the 
chosen representation $\pi\rtimes U$ (e.g., see \cite{CELY}*{Proposition 2.7.4}), 
it suffices to show that there is a unique $x\in X$ such that $\pi(J_x)=\{0\}$ but
$\pi(I_x)\neq \{0\}$. 

As a first step, we want to pass from $A$ to $A/J$ with $J:=\ker\pi$ in order to assume without loss of generality that 
$\pi:A\to \mathcal B(\H_\pi)$ is faithful on $A$. Indeed, since $\Prim(A/J)$ is a closed $G$-invariant subset of $\Prim(A)$, 
and since $\phi:\Prim(A)\to X$ is $G$-invariant, open, and continuous,  one easily checks that $Y:=\phi(\Prim(A/J))$ is closed in $X$ and 
 the restriction $\phi:\Prim(A/J)\to Y$ satisfies all assumptions of the theorem. 
We then have $Q\in \Prim(A/J\rtimes_\alpha G)$ and $(A/J)_x=(I_x+J)/(J_x+J)$ for all $x\in X$. Since $\pi(J_x)=\{0\}$ iff
$\pi(J_x+J)=\{0\}$ and $\pi(I_x)\neq \{0\}$ iff $\pi(I_x+J)\neq\{0\}$, we see that the desired result will pass from $A/J$ to $A$.

So assume now that $\pi$ is faithful.  Then, for every $G$-invariant non-zero ideal  $I\subseteq A$ it follows that
$\pi\rtimes U$ restricts to an irreducible (hence non-degenerate)
representation of the ideal $I\rtimes_\alpha G$ of $A\rtimes_\alpha G$.
Therefore restriction of $\pi$ to $I$ is non-degenerate as well.

Next we show that there exists an element $x\in X$ such that $\{x\}$ is open and dense in $X$. Indeed, since $X$ is almost Hausdorff, it contains a 
dense open Hausdorff subset $U\subseteq X$. If $U$ contains two different points $x_1\neq x_2$, then 
there also exist nonempty disjoint open subsets $U_1, U_2$ of $X$.
Let $I_{U_1}$ and $I_{U_2}$ denote the ideals in $A$ corresponding to the disjoint open sets $\phi^{-1}(U_1)$ and $\phi^{-1}(U_2)$ in $\Prim(A)$,
respectively. It then follows that $I_{U_1}\cdot I_{U_2}=\{0\}$. But since  $I_{U_1}$ and $I_{U_2}$ are both non-zero $G$-invariant ideals,
it follows that
$$\{0\}=\pi(I_{U_1}\cdot I_{U_2})\H_\pi=\pi(I_{U_1})\big(\pi(I_{U_2})\H_\pi\big)=\pi(I_{U_1})\H_\pi=\H_\pi,$$
a contradiction.

So let $x\in X$ be as above. Then $\phi^{-1}(x)$ is dense in $\Prim(A)$ as well and therefore $I_x=\cap\{P: \phi(P)=x\}$ must be $\{0\}$.
On the other hand, if $x\neq y\in X$, then the open set $\phi^{-1}(x)$ lies in the complement of $\phi^{-1}(y)$, and therefore 
$\phi^{-1}(y)$ isn't dense in $\Prim(A)$. Thus $I_y=\ker\big(\phi^{-1}(y))\neq \{0\}$. 
So we see that $x\in X$ is the unique element with $\{0\}=I_x$. 
Since $\pi$ is faithful and $J_x$ is strictly larger than $I_x=\{0\}$, it  also follows that $\pi(J_x)\neq 0$.
\end{proof}

If $A$ is separable and $G$ is second countable, then the following much stronger result holds true:

\begin{proposition}\label{orbit-method-separable}
Suppose that $\alpha:G\car A$ is an action  with $G$ second countable and $A$ separable.
Suppose further that $(A_i)_{i\in I}$ is a decomposition of $A$ by $G$-invariant subquotients.
Then $(A_i\rtimes_{\alpha_i}G)_{i\in I}$ is a decomposition of $A\rtimes_\alpha G$ by subquotients.
\end{proposition}
\begin{proof}
Again, let $Q=\ker(\pi\rtimes U)$ for some  irreducible representation $\pi\rtimes U$ 
of $A\rtimes_{\alpha}G$. By \cite{Green}*{Corollary 19} there exits some $P\in \Prim(A)$ such that 
$\ker\pi=\Res_A(Q)=\ker G(P)$. By assumption, there exists  a unique $i\in I$ such that $P\in \Prim(A_i)$, and 
then $G(P)=\{\alpha_g(P):g\in G\}\subseteq  \Prim(A_i)$ as well.  Let $A_i=I_i/J_i$.
It then follows that that $J_i\subseteq \ker G(P)=\ker\pi$, but 
and $\pi(I_i)\neq \{0\}$ since $I_i$ is strictly larger than $J_i$.  Therefore 
$\pi\rtimes U$ is an irreducible representation of $A_i\rtimes_{\alpha_i}G$.
 \end{proof}
 
 Let us also mention the following (well known) consequence of Propositions \ref{prop-orbit-method} and \ref{orbit-method-separable}

\begin{corollary}\label{cor-orbit}
Suppose that $\alpha:G\car A$ is an action of $G$ on the $C^*$-algebra $A$ such that every orbit $G(P)=\{\alpha_g(P): g\in G\}$ in 
$\Prim(A)$ is locally closed.
Assume further that one of the following assumptions hold
\begin{enumerate}
\item $A$ is separable and $G$ is second countable, or
\item $\Prim(A)/G$ is almost Hausdorff.
\end{enumerate}
Then $(A_{G(P)}\rtimes_{\alpha}G)_{G(P)\in \Prim(A)/G}$ is a decomposition of $A\rtimes_\alpha G$ 
by subquotients. 
%
\end{corollary}

We are coming back to the case of a smooth action $\alpha:G\car A$  in the sense of Definition \ref{def-smooth}.
It follows then from  Corollary \ref{cor-orbit} that 
$A\rtimes_{\alpha}G$ admits  a decomposition by subquotients 
$$(A_{G(P)}\rtimes_{\alpha} G)_{G(P)\in\Prim(A)/G}$$
 over the orbit space $\Prim(A)/G$.
Since we have $G$-equivariant homeomorphisms $\Prim(A_{G(P)})\cong G(P)\cong G/G_P$,
and $G(P)\subseteq \Prim(A)$ is locally closed, it follows in particular that $\{P\}$ is closed in $G(P)$ and hence locally closed in $\Prim(A)$ as well. Let $A_P:=A_{\{P\}}$ denote the corresponding simple subquotient of $A$. 
Since $P$ and $\ker(\overline{\{P\}})$ are invariant under the restriction of $\alpha$ to the stabilizer $G_P$ of $P$, we  obtain a 
canonical action $\alpha^P:G_P\car A_P$ and it
follows  from Theorem \ref{thm-induced} that $\alpha:G\car A_{G(P)}$ is isomorphic to
$\Ind\alpha^P:G\car \Ind_{G_P}^GA_P$. 
By Green's imprimitivity theorem (Theorem \ref{thm-Green}), it follows that 
$A_{G(P)}\rtimes_{\alpha}G$ is Morita equivalent to $A_P\rtimes_{\alpha^P} G_P$.
Moreover, if $G$ is abelian, then it follows from Example \ref{ex-induced-Mor} that this Morita equivalence is equivariant 
with respect to the  dual action 
$\widehat\alpha:\widehat{G}\car A_{G(P)}\rtimes_{\alpha}G$ and the inflation to $\widehat{G}$ of the dual action 
$\widehat{\alpha^P}:\widehat{G_P}=\widehat{G}/G_P^\perp \car A_P\rtimes_{\alpha^P}G_P$.

\section{Twisted group algebras and actions on $\K(\H)$}\label{sec-actionsK}
In this section we want to study smooth actions on $C^*$-algebras of type I. Recall  that if $A$ is of type I, then 
$\widehat{A}\cong \Prim(A)$ via $[\rho]\mapsto \ker\rho$, and then the simple subquotient $A_\rho:=A_{\ker\rho}$ is isomorphic 
to the algebra of compact operators $\K(\H_\rho)$. Thus, by the discussion at the end of the preceding section, we want to understand 
the crossed products $\K(\H_\rho)\rtimes_{\alpha^\rho}G_\rho$ with $G_\rho:=G_{\ker\rho}$ and $\alpha^\rho:G_\rho\car \K(\H_\rho)$ the corresponding action.

Recall  that for any Hilbert space $\mathcal H$ the projective unitary group $P\U(\H):=\U(\H)/\T 1_{\H}$ is isomorphic to
 $\Aut(\K(\H))$ via $[V]\mapsto \Ad V$. Thus if 
$\beta:G\car \K(\H)$ is an action, we can form the group
\begin{equation}\label{eq-Gbeta}
G_\beta:=\{(g, V)\in G\times \U(\H): \beta_g=\Ad V\}.
\end{equation} Then $G_\beta$ fits into a central extension
$$ 1\to \T\to G_\beta\to G\to 1$$
with inclusion map $z\mapsto (e, z1_{\H})$ and quotient map $(g,V)\mapsto g$. 
Choosing then a Borel section $c:G\to G_\beta$ and writing $c(g):=(g, V_g^\beta)$, we obtain a Borel map
$V^\beta:G\to \U(\H)$ and a Borel map $\om_\beta:G\times G\to \T$
 given by
 \begin{equation}\label{eq-omrep}\om_\beta(g,h)1_\H=V^\beta_gV^\beta_h(V^\beta_{gh})^*.
 \end{equation}
It is then easy to check that $\om_\beta$ is a Borel $2$-cocycle on $G$, i.e., it satisfies the cocycle condition
\begin{equation}\label{eq-cocycle}
\om_\beta(g,h)\om_\beta(gh, l)=\om_\beta(g,hl)\om_\beta(h,l)\quad \forall g,h,l \in G.
\end{equation}
Moreover, the isomorphism class $[G_\beta]$ of the group extension 
$1\to \T\to G_\beta\to G\to 1$ and the 
 cohomology class $[\om_\beta]\in H^2(G,\T)$\footnote{Recall that two cocycles $\om,\om'\in Z^2(G,\T)$ are equivalent (or cohomologous), iff there 
 exists a Borel map $f:G\to \T$ such that $\om(g,h)=\partial f(g,h)\om'(g,h)=f(g)f(h)f(gh)^{-1}\om'(g,h)$ for all $g,h\in G$.}  only depend on the Morita equivalence class 
 of the action $\beta:G\car \K(\H)$ (e.g., see \cite{BE:deformation}*{Section 4} for a detailed account in the 
 non-separable case).
 
%
%
%
Conversely,  if $\om\in Z^2(G,\T)$ is a Borel cocycle, then $\om$ determines a  central extension $1\to \T\to G_\om\to G\to 1$ of $G$ by $\T$ 
where $G_\om:= G\times \T$ (as a Borel space) equipped with group operations given by 
 \begin{equation}\label{eq-Gom}
 (g_1, z_1)(g_2, z_2)=(g_1g_2, \om(g_1, g_2)z_1 z_2)\quad\text{and}\quad (g,z)^{-1}=(g^{-1},\overline{\om(g,g^{-1})}\bar z).
 \end{equation}
  It has been shown by Mackey (for second countable $G$) and Kleppner \cite{Klep1} (in the general case) that there is a unique
locally compact group topology on $G_\om$ which induces the product Borel structure\footnote{This does not have to be
the product topology in general}. If $\om_\beta$ is as above,  then  $[G_{\om_\beta}]=[G_\beta]$  as group extensions. 
 
 The inversion formula in $G_\om$ becomes easier if $\om$ is a {\em normalized cocycle} in the sense that 
 \begin{equation}\label{eq-normalized}
 \om(g,g^{-1})=1
 \end{equation}
for all $g\in G$. It then also satisfies the equation
\begin{equation}\label{eq-normalized1}
\om(g, h)=\overline{\om(h^{-1},g^{-1})} 
\end{equation}
for all $g,h\in G$ (see \cite{Klep1}*{p. 216}).
 It is shown in \cite{Klep1}*{p. 215} that every cocycle is equivalent to a {\em normalized cocycle}, so one may always pass to a normalized one, if this seems  convenient. 

There is a canonical Green-twisted action $(\id, \iota^\om)$ of the pair $(G_\om,\T)$ on the complex numbers 
$\C$ in which $\id:G_\om\car \C$ is the trivial action and $\iota^\om:\T\to \C$ is the inclusion map. The twisted group algebra 
$C^*(G,\om)$ can  then be defined as the twisted crossed product $\C\rtimes_{(\id,\iota^\om)}G_\om/\T$. 
The dense subalgebra $C_c(G_\om,\C,\iota^\om)$ as in (\ref{eq-Green-twist2}) then consists of all comapactly supported functions
$f:G_\om\to \C$  satisfying $f(g,z)=\bar{z}f(g,1)$ for all $(g,z)\in G_\om$, equipped with convolution and involution
given by the restriction of the usual convolution and involution on $C_c(G_\om)\subseteq C^*(G_\om)$.
It is  easily checked that $C^*(G,\om)$ depends, up to isomorphism, 
only on the cohomology class of $\om$. Moreover,   there
is a one-to-one correspondence between  
covariant representations for the twisted action $(\id,\iota^\om):(G_\om,\T)\car \C$ and 
 the $\om$-representations $V:G\to\U(\H)$, i.e., Borel maps $V:G\to\U(\H)$ that satisfy 
 \begin{equation}\label{eq-omrepresentation}
 V_gV_h=\om(g,h)V_{gh}\quad \forall g,h\in G.
 \end{equation}
 Indeed, given such $\om$-representation, we obtain a unitary representation $V':G_\om\to \U(\H)$ by putting $V'_{(g, z)}:=zV_g$ that clearly preserves 
 the twist $\iota$ as in (\ref{eq-Green-twist1}). \footnote{Alternatively, $C^*(G,\om)$ can also be defined as the enveloping $C^*$-algebra of the twisted  convolution algebra $L^1(G,\om)$, 
 i.e., $L^1(G)$ equipped with involution and convolution twisted by $\om$ in a certain sense, e.g., see \cite{BK}.}
 
 \begin{remark}\label{rem-quotient}
 It follows from the above constructions that the twisted group algebra $C^*(G,\om)$ is  a direct summand of the algebra 
 $\C\rtimes_{\id} G_\om=C^*(G_\om)$. In particular, it follows that the irreducible representations of twisted group algebras of compact groups are finite dimensional.
 \end{remark}

Now, if $V:G\to \U(\H)$ is an $\om$-representation and if $\beta=\Ad V:G\car \K(\H)$ is the corresponding action on $\K(\H)$, we 
obtain an isomorphism 
$\K(\H)\otimes C^*(G,\bar\om)\cong \K(\H)\rtimes_\beta G$ given as an extension of the mapping
\begin{equation}\label{eq-actionK}
\K(\H)\odot C_c(G_{\bar\om},\C,\iota)\to C_c(G, \K(H)); k\otimes f\mapsto [g\mapsto f(g,e)kV_g^*].
\end{equation}
Here $\bar\om(g,h):=\overline{\om(g,h)}$ denotes the complex conjugate of $\om$.
This is an easy version of  \cite{Green}*{Theorem 18}. Note that equation (\ref{eq-omrep}) above implies that 
$V_\beta:G\car \U(\H)$ is an $\om_\beta$-representation, and therefore we obtain the isomorphism 
$\K(\H)\rtimes_\beta G\cong \K(\H)\otimes C^*(G,\bar\om_\beta)$ for any given action $\beta:G\car \K(\H)$.

\begin{remark}\label{rem-dual-twisted-action}
If $G$ is abelian and $\om\in Z^2(G,\T)$ is a Borel cocycle, we have a dual action, say $\widehat{\om}:=\widehat{(\id,\iota^\om)}:\widehat{G}\car C^*(G,\om)$,
given by 
$$\big(\widehat{\om}_\chi(f)\big)(g, z)=\overline{\chi(g)}f(g,z)\quad\forall \chi\in \widehat{G}, \; (g,z)\in G_\om.$$
This is of course a special case of the dual action as defined in item (3) of 
Remark \ref{rem-twistedMor}. One then easily checks that the isomorphism in (\ref{eq-actionK}) is 
$\id_\K\otimes \widehat{\bar\om}-\widehat{\beta}$ 
equivariant. In particular, we obtain a  $\widehat{\beta}-\widehat{\bar\om_\beta}$ equivariant $\K(\H)\rtimes_\beta G - C^*(G,\bar\om_\beta)$ Morita equivalence
 for every action $\beta:G\car \K(\H)$ of an abelian group $G$.
\end{remark}

Combining the above discussions with Corollary  \ref{cor-orbit}, we arrive at

\begin{theorem}\label{thm-transitve-typeI}
Suppose that $\alpha:G\car A$ is a smooth action of $G$ on the \mbox{type I} $C^*$-algebra $A$.
Then $A\rtimes_{\alpha}G$ admits a decomposition $(A_{G(\rho)}\rtimes_\alpha G)_{G(\rho)\in \widehat{A}/G}$ 
such that each subquotient $A_{G(\rho)}\rtimes_\alpha G$ is Morita equivalent to a twisted group algebra $C^*(G_\rho,\om_\rho)$.
Moreover, if $G$ is abelian, the dual actions $\widehat{\alpha}:\widehat{G}\car A_{G(\rho)}\rtimes_\alpha G$
are Morita equivalent to the inflations to $\widehat{G}$ of the dual actions 
$\widehat{\om_\rho}:\widehat{G_\rho}\car C^*(G_\rho, \om_\rho)$.

\end{theorem}

By applying Theorem \ref{thm-Morita} we obtain a similar result for smooth twisted actions $(\alpha,\tau):(G,N)\car A$ on  type I $C^*$-algebras $A$. For later use, we close this section with the following well-known result:

\begin{theorem}\label{thm-compact-typyI}
Suppose that $\alpha:G\car A$ is an action of the compact group $G$ on the type I $C^*$-algebra $A$. Then $A\rtimes_\alpha G$ is  type I as well.
\end{theorem}
\begin{proof} It follows from  \cite{Ungermann}*{Theorem 27} that $\alpha:G\car A$ is smooth. The result then follows from 
Theorem \ref{thm-transitve-typeI} and the fact that all twisted group algebras $C^*(G_\rho,\om_\rho)$ for the compact stabilizers $G_\rho$ for the action $G\car \widehat{A}$ are type I (see Remark \ref{rem-quotient}).
\end{proof}
%
%
%
%
%

%
%
%

\section{The main results}\label{sec-connected}

We are now going to prove our main results, which will imply that the simple subquotients of certain crossed products by actions of  abelian groups  are Morita equivalent to twisted group algebras of  abelian groups.
We start with the following observation:

 \begin{proposition}\label{prop-dual-smooth}
 Suppose that $\alpha:G\car A$ is a smooth action of the abelian locally compact group $G$ on the type I  $C^*$-algebra $A$. 
  Then the orbit spaces $\Prim(A\rtimes_{\alpha}G)/\widehat{G}$ and $\Prim(A)/G\cong \widehat{A}/G$ are homeomorphic and the dual action $\widehat{\alpha}:\widehat{G}\car A\rtimes_{\alpha}G$ is smooth as well.
  \end{proposition}
  \begin{proof} We need to show that the induced action $\widehat{G}\car \Prim(A\rtimes_{\alpha}G)$ satisfies the requirements of
  Definition \ref{def-smooth}. Indeed, by Theorem \ref{thm-transitve-typeI} we see that $A\rtimes_\alpha G$ admits a decomposition 
  by subquotients $(A_{G(\rho)}\rtimes_\alpha G)_{G(\rho)\in \widehat{A}/G}$ such that each fibre 
  $A_{G(\rho)}\rtimes_\alpha G$ is $\widehat{G}$-equivariantly Morita equivalent to some twisted group algebra 
  $C^*(G_\rho,\om_\rho)$ for the stabilizer $G_\rho$ of $\rho$ and some cocycle $\om_\rho\in Z^2(G_\rho,\T)$, where 
  the $\widehat{G}$-action on $C^*(G_\rho,\om_\rho)$ is the inflation of the dual action $\widehat{\om_\rho}:\widehat{G_\rho}\car C^*(G_\rho,\om_\rho)$.  It follows then from Remark \ref{rem-action-Prim} that the action 
  $\widehat{G}\car \Prim(C^*(G_\rho,\om_\rho))\cong \widehat{S_{\om_\rho}}$ is transitive with common stabilizer 
  $\widehat{G}_P=S_{\om_\rho}^{\perp}\subseteq \widehat{G}$. 
  So we see that the $\widehat{G}$-orbits $\widehat{G}(P)$ in $\Prim(A\rtimes_{\alpha}G)$ coincide with the locally closed subsets 
  $\Prim(A_{G(\rho)}\rtimes_{\alpha}G)$ which are equivariantly  homeomorphic  to  
  $ \widehat{S_{\om_\rho}}\cong \widehat{G}/S_{\om_\rho}^\perp\cong \widehat{G}/\widehat{G}_P$.

In particular, from the above we obtain a bijection $\widehat{A}/G\longleftrightarrow \Prim(A\rtimes_{\alpha}G)/\widehat{G}$ sending 
an orbit $G(\rho)\in \widehat{A}/G$ to the orbit $\Prim(A_{G(\rho)}\rtimes_\alpha G)\in \Prim(A\rtimes_\alpha G)/G$. Using continuity of induction and restriction of representations (e.g., see \cite{CELY}*{Proposition 2.7.4}), one even checks that this bijection is continuous in both directions (in the separable case, this follows from \cite{GootLaz}). This finishes the proof.  \end{proof}

\begin{remark}\label{rem-dualsmooth}
Using  item (3) of  Remark \ref{rem-twistedMor}, we can easily swap the roles of $\alpha$ and $\widehat{\alpha}$ in the above proposition: so if $A\rtimes_{\alpha}G$ is type I and $\widehat{\alpha}:\widehat{G}\car A\rtimes_\alpha G$ is smooth, then 
$\alpha:G\car A$ is  smooth and the respective orbit spaces are homeomorphic.
\end{remark}

Before we state our main theorem, we also need to observe the following

\begin{lemma}\label{lem-twisted-crossed}
Let $\om$ be a normalized Borel $2$-cocycle on the abelian group $L$ and let $\varphi:G\to \widehat{L}$ be a continuous group homomorphism.
Let $\beta:= \widehat{\om}\circ \varphi: G\car C^*(L,\om)$ be the composition of the dual action $\widehat{\om}:\widehat{L}\car C^*(L,\om)$   with $\varphi$. Then
$$C^*(L,\om)\rtimes_\beta G\cong C^*(G\times L, \eta)$$
for the normalized cocycle $\eta\in Z^2(G\times L,\T)$ given by 
$$\eta\big((g, l), (g', l')\big)=\varphi(g)(l')\om(l,l')\quad \forall (g,l), (g',l')\in G\times L.$$
\end{lemma}
\begin{proof}
It is easy to check that $\eta$ is a Borel $2$-cocycle on $G\times L$ and that 
$$\Phi: C_c(G, C_c(L_\om, \C, \iota^\om))\to C_c((G\times L)_\eta, \C, \iota^\eta)\; 
\quad \Phi(f)\big( ((g,l), z)\big)=f(g)(l,z)$$
extends to the desired isomorphism. 
\end{proof}
 
  Our central result is now the following:

\begin{theorem}\label{thm-main}
Suppose that $\gamma:L\car A$ is an action of the abelian  locally compact group $L$ on the $C^*$-algebra $A$ such that 
$A\rtimes_\gamma L$ is type I and the dual action $\widehat\gamma:\widehat{L}\car A\rtimes_\gamma L$ is smooth.
Let $\psi:G\to L$ be a continuous 
homomorphism from the locally compact group $G$ to $L$, and let  $\alpha=\gamma\circ \psi:G\car A$. 

Then there is a decomposition by subquotients $(A_{L(P)}\rtimes_\alpha G)_{L(P)\in \Prim(A)/L}$ of 
$A\rtimes_\alpha G$ such that each fibre $(A\rtimes_\alpha G)_{L(P)}$ 
is Morita equivalent to some twisted group algebra $C^*(G\times H_\rho, \eta_\rho)$, where $H_\rho\subseteq \widehat{L}$ 
is the stabilizer of some $\rho\in (A_{L(P)}\rtimes_\beta L)\dach$ under the dual action $\widehat{\beta}$.

\end{theorem}
\begin{proof} 
It follows from Remark \ref{rem-dualsmooth} that the action $L\car \Prim(A)$ is smooth, hence $X:=\Prim(A)/L$ is almost Hausdorff and the orbits $L(P)$ are locally closed in $\Prim(A)$. Moreover, it follows from the construction of $\alpha=\gamma\circ \psi$ that 
the subquotients $A_{L(P)}$ are $G$-invariant. Thus, 
it follows from Proposition \ref{prop-orbit-method} that $A\rtimes_\alpha G$ has a decomposition by subquotients
$(A_{L(P)}\rtimes_\alpha G)_{L(P)}$ over the orbit space $\Prim(A)/L$.
So we only need to show that every subquotient $A_{L(P)}\rtimes_\alpha G$ is Morita equivalent to some twisted group algebra 
$C^*(G\times H_\rho,\om_\rho)$ as in the theorem.

Since the statement of the theorem is invariant under passing to Morita equivalent actions, we may pass from 
$\gamma$ (resp. $\alpha=\gamma\circ \psi$) to $\widehat{\widehat{\gamma}}$ (resp. $\alpha':= \widehat{\widehat{\gamma}}\circ \psi$). 
In particular, we obtain a Morita equivalence
$$A_{L(P)}\rtimes_\alpha G\sim_M \big((A_{L(P)}\rtimes_\gamma L)\rtimes_{\widehat{\gamma}}\widehat{L}\big)\rtimes_{\alpha'}G.$$

By our assumptions and by Proposition \ref{prop-dual-smooth} and Remark \ref{rem-dualsmooth} we know that $A_{L(P)}\rtimes_\gamma L$ is of type I and that $(A_{L(P)}\rtimes_\gamma L)\dach$ is canonically homeomorphic to a homogeneous 
$\widehat{L}$-space $\widehat{L}/H_\rho$, where $H_{\rho}\subseteq \widehat{L}$ denotes the stabilizer of some  
$\rho\in(A_{L(P)}\rtimes_\gamma L)\dach$. By Theorem \ref{thm-transitve-typeI} we see that
$(A_{L(P)}\rtimes_\gamma L)\rtimes_{\widehat{\gamma}}\widehat{L}$ is $L$-equivariantly Morita equivalent to a twisted group algebra 
$C^*(H_\rho, \om_\rho)$, where the $L$-action 
on $C^*(\widehat{L}_\rho,\om_{\rho})$ is via the inflation of the dual 
action $\widehat{\om_\rho}: \widehat{H_\rho}\cong L/H_{\rho}^\perp \car C^*(\widehat{L}_\rho,\om_{\rho})$ to $L$. 
It then follows that $\alpha:G\car A_{L(P)}$ is Morita equivalent to the action $\beta: G\car C^*(H_\rho,\om_{\rho})$
given by the composition $\beta=\widehat{\om_\rho}\circ \varphi$, where $\varphi:G\to \widehat{H_\rho}$ 
is the composition of group homomorphisms 
$$G\stackrel{\psi}{\longrightarrow} L\stackrel{q}{\longrightarrow} L/H_\rho^\perp\cong \widehat{H_\rho}.$$
It follows then from Lemma \ref{lem-twisted-crossed} that 
$A_{L(P)}\rtimes_{\alpha}G$ is Morita equivalent to 
$C^*(H_\rho,\om_\rho)\rtimes_\beta G\cong C^*(G\times H_\rho,\eta_\rho)$ with $\eta_\rho$ given by
$$\eta_\rho((g, \chi), (g,\chi'))=\chi'(\varphi(g))\om_\rho(\chi,\chi')$$
for all $(g,\chi), (g',\chi')\in G\times H_\rho$.
\end{proof}

Using the Morita equivalence between $\gamma$ and the double dual action $\widehat{\widehat{\gamma}}$, we can also formulate the following dual version of Theorem \ref{thm-main}. 

\begin{theorem}\label{thm-main1}
Suppose that $\gamma:L\car A$ is a smooth action of the abelian locally compact group $L$ on the type I $C^*$-algebra $A$ and let 
$\psi: G\to \widehat{L}$ be a continuous homomorphism from the locally compact group $G$ to $\widehat{L}$. 
Let 
$$\beta:=\widehat{\gamma}\circ \psi: G\car B:=A\rtimes_\beta L.$$
Then $B\rtimes_\beta G$ has a decomposition 
by subquotients $(B_{\widehat{L}(P)}\rtimes_\beta G)_{\widehat{L}(P)\in \Prim(B)/\widehat{L}}$ such that for each fibre 
$B_{\widehat{L}(P)}\rtimes_\beta G$ the following holds:
\begin{enumerate}
\item There exists a unique $L$-orbit $L(\rho)$ in $\widehat{A}$ such that $B_{\widehat{L}(P)} \cong A_{L(\rho)}\rtimes_\gamma L$, and 
\item $B_{\widehat{L}(P)}\rtimes_\beta G$ is Morita equivalent to some twisted group algebra $C^*(G\times L_\rho,\om_\rho)$,
where $L_\rho$ denotes the stabilizer of $\rho$ under the action $L\car \widehat{A}$.
\end{enumerate}
\end{theorem}

\begin{remark}\label{rem-smooth} It follows directly from  Proposition \ref{prop-dual-smooth} that if $\gamma:L\car A$ is an action of the abelian locally compact group $L$ on the $C^*$-algebra $A$ such that $A$ and $A\rtimes_\gamma L$ are both of type I, then $\gamma:L\car A$ is smooth if and only if $\widehat{\gamma}:\widehat{L}\car A\rtimes_\gamma L$ is smooth. Hence in this situation, both 
versions, Theorem \ref{thm-main} and Theorem \ref{thm-main1} will apply.
\end{remark}

As a first application of Theorem \ref{thm-main} we now look at the case of connected and simply connected Lie groups as discussed in the introduction:

\begin{corollary}\label{cor-simplyconnected}
Suppose that $G$ is a connected and simply conectetd Lie group. Then $C^*(G)$ admits a decomposition $(B_\lambda)_{\lambda\in \Lambda}$ by subquotients $B_\lambda$ such that each fibre $B_\lambda$ is Morita equivalent to some a twisted group algebra 
$C^*(G_\lambda, \om_\lambda)$, where each $G_\lambda$ is a compactly generated abelian group.
\end{corollary}
\begin{proof}
Following the discussion in the introduction, we find an inclusion $G\into L$ for some simply connected Lie group $L$,
such that $N:=[G,G]=[L,L]$, $C^*(N)$ and $C^*(L)$ are type I, and the conjugation action $\beta:L\car \widehat{N}\cong \widehat{C^*(N)}$ is smooth. Then, if $(\gamma,\tau):(L,N)\car C^*(N)$ denotes the conjugation action, and if $(\alpha,\tau):(G,N)\car C^*(N)$ denotes its restriction to $(G,N)$, we get 
$$C^*(L)\cong C^*(N)\rtimes_{(\gamma,\tau)}L/N\quad\text{and}\quad C^*(G)\cong C^*(N)\rtimes_{(\alpha,\tau)}G/N.$$
Now, using Theorem \ref{thm-Morita} we may pass to a Morita equivalent action $\tilde\gamma:L/N\car A$ and its restriction 
$\tilde\alpha:G/N\car A$ (which is just the composition $\gamma\circ \varphi$ with $\varphi:G/N\into L/N$ the inclusion map).
 Then $A\sim_M C^*(N)$ and $A\rtimes_{\tilde\gamma}L/N\sim_M C^*(L)$ are type I and  $\tilde\gamma:L/N\car A$ is smooth
(since $L/N\car \widehat{A}$ is equivariantly homeomorphic to $L/N\car \widehat{N}$). By Remark \ref{rem-smooth} it follows that
that  Theorem \ref{thm-main} applies to  $A\rtimes_{\tilde\alpha}G/N$. It follows that  $A\rtimes_{\tilde\alpha}G/N$ 
has a decomposition by subquotients
$(B_\lambda)_{\lambda\in \widehat{A}/L}$ 
such that each $B_\lambda$ is Morita equivalent to a twisted group algebra $C^*(G/N\times H_\rho,\om_\rho)$, where $
H_\rho\subseteq\widehat{L/N}$ is the stabilizer of some $\rho\in (A\rtimes_{\tilde\gamma}L/N)\dach$ under the dual action.
Since $G/N$  and $L/N$ are  connected abelian Lie groups,  it follows 
that $G_\lambda:=G/N\times H_\rho$ is compactly generated. Since $C^*(G)$ is Morita equivalent to $A\rtimes_{\tilde\alpha}G/N$,
the result follows.
\end{proof}

Another nice consequence of Remark \ref{rem-smooth} together with Theorem \ref{thm-compact-typyI} is the following

\begin{corollary}\label{cor-compact}
Suppose that $\gamma:K\car A$ is an action of the compact group $K$ on the type I $C^*$-algebra $A$. 
Suppose further that $G$ is an abelian locally compact group and that $\psi:G\car K$ is a continuous 
group homomorphism. Then Theorem \ref{thm-main} applies to the crossed product $A\rtimes_\alpha G$
for $\alpha=\gamma\circ \psi$.
\end{corollary}

%
%
%

\section{Twisted group algebras of  abelian groups}\label{sec-twisted-group-algebras}

In this section we want to recall  the structure of twisted group algebras  of abelian locally compact groups $G$. 
We start with  the following  useful (and certainly well-known)  fact for equivalences of cocycles related to closed subgroups.
In what follows, if $\om\in Z^2(G,\T)$ and if $H$ is a closed subgroup of $G$, then we denote by $\om_H$ the restriction 
$\om|_{H\times H}$ of $\om$ to $ H\times H$. We  call it the {\em restriction of $\om$ to $H$}.

\begin{proposition}\label{prop-cocycle}
Let $\om$ be a  Borel $2$-cocycle on the locally compact group $G$ and let $H$ be a closed subgroup of $G$.
Suppose that  $\om_H$ is equivalent to the cocycle $\eta\in Z^2(H,\T)$.
Then there exists a cocycle $\tilde\om\in Z^2(G,\T)$ equivalent to $\om$  such that $\tilde\om_H=\eta$.
If $\om$ and $\eta$ are normalized, then $\tilde\om$ can be chosen to be normalized as well.
\end{proposition}
\begin{proof}
Let $f:H\to \T$ be a Borel map such that $\partial f\cdot \om_H=\eta$. Let $\tilde{f}:G\to \T$ be defined  
by $\tilde{f}|_H=f$ and $\tilde{f}(g)=1$ for $g\notin H$. Then $\tilde\om :=\partial\tilde{f}\cdot \om$ satisfies $\tilde\om_H=\eta$.
For the final statement simply observe that if $\om$ and $\eta$ satisfy (\ref{eq-normalized}), then so does $\tilde\om$. 
\end{proof}
%
%
%
%
%
%
%
%

We proceed by recalling some facts on cocycles of abelian groups due to Kleppner \cite{Kleppner} and Baggett and Kleppner \cite{BagKlep}.
For this let $\om\in Z^2(G,\T)$ with $G$ abelian.
The cocycle identities (\ref{eq-cocycle}) then  imply 
 that  $\om$ induces a homomorphism
\begin{equation}\label{eq-hom}
h_{\om}:G\to \widehat{G}; \;h_{\om}(g)(h):=\om(g,h)\overline{\om(h,g)}
\end{equation}
as studied by Kleppner in \cite{Kleppner}*{Section 7}. A cocycle $\om$ is called {\em totally skew} (or non-degenerate)
if the group kernel $S_\om:=\ker h_\om$ is trivial. Note that 
$$S_\om=\{s\in G: \om(s,g)=\om(g,s)\;\forall g\in G\}$$
is called the  {\em symmetry group} for $\om$. If 
$\om_1\sim\om_2$, then $h_{\om_1}=h_{\om_2}$, and hence $S_{\om_1}=S_{\om_2}$. 
The following result is due to  Baggett and Kleppner
%
%

\begin{theorem}[\cite{BagKlep}*{Theorem 3.1}]\label{thm-lift}
Let $G$ be an abelian locally compact group and let  $\om \in Z^2(G,\T)$. Then  there exists a normalized totally skew Borel cocycle $\tilde{\om}\in Z^2(G/S_\om,\T)$ 
such that $\om$ is equivalent to the inflated cocycle $\tilde{\om}\circ (q\times q)\in Z^2(G,\T)$, where 
$q:G\to G/S_\om$ denotes the quotient map. 
\end{theorem}

In particular, we may assume without loss of generality (by passing  to an equivalent cocycle if necessary), that 
\begin{equation}\label{eq-omsym}
\om(g,s)=1=\om(s, g)\quad \forall g\in G, s\in S_{\om}.
\end{equation}
We now want to consider the dual action $\widehat{\om}:\widehat{G}\car C^*(G,\om)$ given by $\widehat{\om}_\chi(f)=\bar{\chi}\cdot f$  for $f\in C_c(G_\om,\C,\iota^\om)$.
By the above theorem we may assume that $\om$ is inflated from a normalized cocycle $\tilde{\om}\in Z^2(G/S_\om,\T)$ as above. 
We want to prove

\begin{theorem}\label{thm-induced-twisted-group-algebra}
Let $\om\in Z^2(G,\T)$ be as above. Let $S:=S_\om$. Then there is an isomorphism of dynamical systems
$$(C^*(G,\om), \widehat{G}, \widehat{\om})\cong \left(\Ind_{\widehat{G/S}}^{\widehat{G}}C^*(G/S,  \tilde\om), \widehat{G}, \Ind\widehat{\tilde\om}\right).$$
\end{theorem}
\begin{proof} Let $S^\perp:=\{\chi\in \widehat{G}: \chi|_S\equiv 1\}$. We want to apply Theorem \ref{thm-induced} above. For this we first construct a 
$\widehat{G}$-equivariant continuous map
\begin{equation}\label{eq-restriction-action}
\varphi:C^*(G,\om)\dach\to \widehat{S}=\widehat{G}/S^\perp. 
\end{equation}
Recall that the irreducible representations of $C^*(G,\om)$ are just the integrated forms of the 
irreducible $\om$-representations $V:G\to \U(\H_V)$. 
We claim that for all such $V$ there exists a unique character $\chi_V\in \widehat{S}$ such that 
$V|_S=\chi_V\cdot 1_{\H_V}$. Indeed, this follows from Schur's lemma, since for all $s\in S, g\in G$
we get
$$V_s V_g=\om(s,g)V_{sg}\stackrel{(\ref{eq-omsym})}{=}V_{sg}=V_{gs}\stackrel{(\ref{eq-omsym})}{=}\om(g,s)V_{gs}=V_gV_s,$$
and hence Schur's lemma implies that $V_s=\chi_V(s)1_{\H_V}$ for some element $\chi_V(s)\in \T$.
Since $V|_S$ is a (Borel)-homomorphism (use that $\om|_{S\times S}\equiv 1$) it follows that 
$\chi_V\in \widehat{S}$. Since restriction of representations is continuous it follows that $\varphi: V\mapsto \chi_V$ is continuous.

To see that it is $\widehat{G}$-equivariant, first observe that the action $\widehat{G}\car  C^*(G,\om)\dach$ 
induced from $\widehat{\om}:\widehat{G}\car C^*(G,\om)$ is given by
$(\nu, V)\mapsto \nu\cdot V$
(pointwise multiplication), since for all $f\in C_c(G_\om, \C,\iota^\om)$ we have
$$V(\widehat{\om}_{\nu^{-1}}(f))=\int_G \widehat{\om}_{\nu^{-1}}(f)(g)V_g\,dg=\int_G \nu(g)f(g)V_g\, dg=(\nu \cdot V)(f).$$
Hence we get  $\varphi(\nu V)=\chi_{(\nu V)}=\nu\chi_V=\nu \varphi(V)$ for $\nu\in \widehat{G}$, $V\in C^*(G,\om)\dach$.

We now look at those $V\in C^*(G,\om)\dach$ which map to $1_S$. Then $V|_S$ is trivial and $V$ is  inflated  from a $\tilde{\om}$-representation $\tilde{V}:G/S\to \U(\H_V)$
defined by $\tilde{V}_{gS}:=V_g$. Note that $\tilde V$  is well defined since 
for all $g\in G$ and $s\in S$ we have 
$$V_{gs}=\bar{\om}(g,s)V_gV_s=\bar{\om}(g,s)V_g\stackrel{(\ref{eq-omsym})}{=}V_g.$$
It follows  that the ideal
$$I=\cap\{\ker V: \varphi(V)=1_S\}$$
is the kernel of the quotient map $Q:C^*(G,\om)\to C^*(G/S, \tilde{\om})$ defined on the level of 
functions by sending $f\in C_c(G_\om,\C, \iota^\om)$ to $Q(f)\in C_c((G/S)_{\tilde\om},\C, \iota^{\tilde\om})$ given by 
  $Q(f)(gS, z)=\int_S f(gs,z)\, ds$. 
The ideal $I$ is invariant under the restriction of $\widehat\om$ to $S^\perp= \widehat{G/S}$ and the 
isomorphism $C^*(G,\om)/I\cong C^*(G/S,\tilde\om)$ induced by $Q$ identifies $\widehat{\om}|_{\widehat{G/S}}$ with $\widehat{\tilde\om}$.
The result now follows from  Theorem \ref{thm-induced}.
\end{proof}

\begin{remark}\label{rem-twisted-group}
It follows from the above theorem, that the twisted group algebra $C^*(G,\om)$ for any abelian locally compact group $G$
is a continuous field of $C^*$-algebras over $\widehat{S}$, for $S:=S_\om$, with fibres isomorphic to the twisted group algebra 
$C^*(G/S,\tilde\om)$ for a fixed totally skew cocycle $\tilde\om\in Z^2(G/S, \T)$ which inflates to a cocycle similar to $\om$ on $G$. 
\end{remark}

The following result is basically a consequence of Kaniuth's \cite{Kaniuth}*{Lemma 2}. Combined with Theorem \ref{thm-induced-twisted-group-algebra}
above it implies the well-known fact that for any cocycle $\om\in Z^2(G,\T)$ on a locally compact abelian group $G$ we have
$\Prim(C^*(G,\om))\cong\widehat{S_\om}$ is Hausdorff.

\begin{theorem}[Kaniuth]\label{thm-Kaniuth}
Let $\om$ be a totally skew Borel 2-cocycle on the locally compact abelian group $G$. 
Then $C^*(G,\om)$ is simple.
\end{theorem}
\begin{proof} 
By the construction of $C^*(G,\om)=\C\rtimes_{(\id, \iota^\om)}G_\om/\T$, we may realize $C^*(G,\om)$ as the quotient of $C^*(G_\om)=\C\rtimes_{\id}G_\om$ by the ideal
$$I_\om:=\cap\{\ker V: V\; \text{a unitary rep. of}\; G_\om, V|_\T=\id_\T\}.$$
Since $\om$ is totally skew, the center $Z(G_\om)$ of $G_\om$ equals $\T$ and if we  apply \cite{Kaniuth}*{Lemma 2} to the 
character $\lambda=\id_\T$, we see that for all $V\in \widehat{G_\om}$ with $V|_\T=\id_\T\cdot 1$ we have
$\ker V=\ker(\Ind_\T^{G_\om} \lambda$). It follows that the quotient $C^*(G,\om)$ has only one 
primitive ideal, and hence it must be simple.
\end{proof}

\begin{remark}\label{rem-action-Prim} In the proof of Theorem  \ref{thm-induced-twisted-group-algebra} we observed 
that the action $\widehat{\om}:\widehat{G}\car C^*(G,\om)\dach$ corresponding to the  dual action $\widehat{\om}$ is given by $(\chi,V)\mapsto \chi\cdot V$ (point-wise multiplication).
Theorem \ref{thm-Kaniuth} now implies that the restriction map $\varphi: C^*(G,\om)\dach\to \widehat{S_\om}$ 
of (\ref{eq-restriction-action}) induces a homeomorphism $\tilde\varphi:\Prim(C^*(G,\om))\to \widehat{S_\om}$.
So, via this isomorphism, we can identify the action $\widehat{G}\car \Prim(C^*(G,\om))$ induced by the dual action 
with the action $\widehat{G}\car \widehat{S_\om}, (\chi,\mu)\mapsto \chi|_{S_\om}\cdot \mu$.
This observation will play an important r\^ole later on.
\end{remark}

In what follows next, we want to study a bit closer the structure of $C^*(G,\om)$ for a totally skew 
cocycle $\om\in Z^2(G,\T)$. In particular, we shall present a proof of Poguntke's theorem
\cite{Pogi}*{Theorem 1} which is equivalent to the result that for compactly generated $G$, the twisted group algebra 
$C^*(G,\om)$ by a totally skew cocycle will be either Morita equivalent to $\C$ or  to
some simple non-commutative torus $C^*(\Z^l,\om_\om')$, of some dimension  $l\geq 2$.
Indeed, we shall see that in the latter case there exists a copy of a finitely generated free abelian group $Z$ inside 
$G$, such that $C^*(G,\om)$ is isomorphic to $C^*(Z,\om_Z)\otimes \K(\H)$ for some separable Hilbert space $\H$.
Although most central ideas are taken from Poguntke's original proof, we hope that our presentation 
streamlines some of the arguments. 

We start with the following observation: suppose that $\om$ is a  Borel 2-cocycle on the abelian locally compact group $G$
and let $H\subseteq G$ be a closed subgroup of $G$.
Then, using Green's decomposition (see Example \ref{ex-decomp}) for twisted crossed products, we can write
$C^*(G,\om)$ as an iterated crossed product 
\begin{equation*}
\begin{split}
C^*(G,\om)&\cong \C\rtimes_{(\id,\iota^\om)}G_\om/\T\cong \big(\C\rtimes_{(\id,\iota^\om)}H_\om/\T\big)\rtimes_{(\alpha,\tau)}G_\om/H_\om\\
&\cong C^*(H,\om_H)\rtimes_{(\alpha,\tau)}G_\om/H_\om,
\end{split}
\end{equation*}
where $H_\om\subseteq G_\om$ denote the central extensions associated to $\om$.
For later use, it is important for us to understand the  twisted action $(\alpha,\tau):(G_\om, H_\om)\car C^*(H,\om_H)$. The following lemma can be seen as an analogue of \cite{Pogi}*{Lemma 1}.

\begin{lemma}\label{lem-action}
In the above situation, if $\om$ is normalized, the decomposition action $\alpha: G_\om\car C^*(H)$ is given by the formula
$\alpha:=\widehat{\om_H}\circ \Res_H^G\circ h_\om^{-1}$, 
where $\widehat{\om_H}:\widehat{H}\car C^*(H,\om_H)$ denotes the dual action, $\Res_H^G:\widehat{G}\to \widehat{H}$  denotes the restriction homomorpphism,
and $h_\om:G_\om\to \widehat{G}$ is (by abuse of notation) the inflation of $h_\om:G\to \widehat{G}$ to $G_\om$.
 The  twist $\tau:H_\om\to U\M(C^*(H,\om_H))$
is given for $f\in C_c(H_\om,\C,\tau)$ by the formula
$$\big(\tau_{(h,z)}(f)\big)(l,v)=f(h^{-1}l, \bar{z}v)\quad \forall  (h,z), (l,v)\in H_\om.$$
\end{lemma}
\begin{proof}
 By Example \ref{ex-decomp} we know that 
  $\alpha:G_\om \car \C\rtimes_{(\id,\iota^\om)}H_\om/\T$ is given for $(g,z)\in G_\om$ and $f\in C_c(H_\om,\C,\iota^\om)$ by the formula
$$\big(\alpha_{(g,z)}(f)\big)(h,v)=f((g,z)^{-1}(h,v)(g,z)),\quad \forall (h,v)\in H_\om.$$
Using (\ref{eq-normalized}) and (\ref{eq-normalized1})
and the fact that $G$ is abelian, we compute
\begin{align*}
(g,z)^{-1}(h,v)(g,z)(h,v)^{-1}&=(g^{-1},\bar{z})(h,v)(g, z)(h^{-1},\bar{v})\\
&=\big(g^{-1}h, \om(g^{-1},h)\bar{z}v\big)\big(gh^{-1},\om(g,h^{-1})z\bar{v}\big)\\
&=\big(g^{-1}hgh^{-1}, \om(g^{-1}h, gh^{-1})\om(g^{-1},h)\om(g,h^{-1})\big)\\
&=\big(e, \om(g^{-1},h)\overline{\om(h, g^{-1})}\big)\\
&=\big(e, h_\om(g^{-1})(h)\big)
\end{align*}
which then implies 
\begin{align*}
\big(\alpha_{(g,z)}(f)\big)(h,v)&=f((g,z)^{-1}(h,v)(g,z))=f(h, [h_\om(g^{-1})(h)] v)\\
&=\big[h_\om(g)(h)\big]f(h,v) =\big(\widehat{\om_H}_{h_\om^{-1}(g)}f\big)(h,v).
\end{align*}
This gives the desired description for the action. 
The formula for the twist follows easily from the general formula as given in 
Example \ref{ex-decomp}.
\end{proof}

\begin{remark}\label{rem-omtrivial}
If in the above lemma the cocycle $\om_H$ is equal to the  trivial cocycle $1_{H\times H}$, then 
$C^*(H,\om_H)\cong C_0(\widehat{H})$, where the  isomorphism extends the Fourier $*$-homomorphism 
$$\Phi:C_c(H_\om, \C, \iota^\om)\to C_0(\widehat{H}); \Phi(f)(\chi)=\int_H f(h,1)\chi(h)\, dh.$$ 
It is easy to check that the decomposition action of Lemma \ref{lem-action} transforms under this isomorphism to 
$(\alpha,\tau):(G_\om,H_\om)\car C_0(\widehat{H})$, where $\alpha$ is induced  by the action
$G_\om\car \widehat{H}$ given via the homomorphism $\Res_H^G\circ h_\om:G_\om\to \widehat{H}$  and 
$$
\big(\tau_{(h,z)} f)(\chi)
= \big(z\chi(h)\big){f}(\chi), \quad f\in C_0(\widehat{H}).
$$
\end{remark}
\medskip

For the next lemma, we need a little preparation about twisted group algebras of inflated cocycles. Assume that 
$H$ is a closed subgroup of the abelian group $G$ and $\om\in Z^2(G,\T)$ is equal to the inflation
of some 
normalized cocycle $\om'\in Z^2(G/H,\T)$. We claim that in this case we get an isomorphism 
$$C^*(G/H,\om')=\C\rtimes_{(\id, \iota^{\om'})}(G/H)_{\om'}/\T\cong  \C\rtimes_{(\id,\sigma)} G_\om/H_\om,$$
where $\sigma:H_\om\to \T=U(\C)$ is given by $\sigma_{(h,z)}=z$.
The first equality is just by definition. To explain 
 the second isomorphism we first observe that $H_\om=H\times \T$ as a group, since $\om_H=1_{H\times H}$ is trivial.
Identifying $H$ with $H\times\{1\}\subseteq G_\om$, we obtain an isomorphism 
$\psi:G_\om/H\to (G/H)_{\om'}; \psi((g,z)H)=(gH,z)$, using that $\om$ is inflated from $\om'$.
Finally, since $\sigma$ is trivial on $H$, it follows that the twisted action $(\id, \sigma)$ of the pair $(G_\om,H_\om)$ is 
inflated from the twisted action $(\id, \iota^{\om'})$ of the pair $((G/H)_{\om'}, \T)$.  The second isomorphism then follows 
by a straightforward  generalization of Example \ref{ex-inflate}.

\begin{lemma}\label{lem-action-omtrivial}
Suppose that $\om$ is a normalized Borel $2$-cocycle on the abelian locally compact group $G$ and that $H\subseteq G$ is a closed subgroup such that the restriction $\om_H$ of $\om$ to $H$ is trivial. Assume further that $\Res_H^G\circ h_\om:G\to \widehat{H}$ is surjective and let $G_H:=\ker(\Res_H^G\circ h_\om)\subseteq G$ and $G'=G_H/H$. Then the restriction of $\om$ to $G_H$ is equivalent to the inflation  of  some  cocycle $\om'\in Z^2(G',\T)$ and then $C^*(G,\om)$ is Morita equivalent  to 
$C^*(G',\om')$.

 If, in addition,  there exists a measurable cross section $G/G_H\to G$ (in particular, if $G$ is second countable), then $C^*(G,\om)$ is isomorphic to $C^*(G',\om')\otimes \K(L^2(\widehat{H}))$.
\end{lemma}
\begin{proof}
If  $G_H\subseteq G$ is as in the lemma, let us write $\tilde{\om}=\om_{G_H}$ and let $h_{\tilde{\om}}:G_H\to \widehat{G_H}$ 
be the corresponding homomorphism. Note that for all $g\in G_H$ we have $h_{\tilde{\om}}(g)=h_\om(g)|_{G_H}$. Hence it follows from the definition of $G_H$ 
that $H\subseteq S_{\tilde\om}\subseteq G_H$. It follows then from Theorem \ref{thm-lift} that there exists a normalized cocycle 
$\om'$ on $G'=G_H/H$ such that $\tilde{\om}$ is equivalent to the inflation, say $\inf\om'$, of $\om'$  to $G_H$. By Proposition \ref{prop-cocycle} we may assume, modulo replacing $\om$ by an equivalent cocycle,  that  $\tilde\om=\inf\om'$.
We now write $C^*(G,\om)$ as iterated twisted crossed 
$C_0(\widehat{H})\rtimes_{(\alpha,\tau)}G_\om/H_\om$ with twisted action $(\alpha,\tau)$ as in as in Remark \ref{rem-omtrivial}.
Since  $G_\om/(G_H)_\om$ is equivariantly homeomorphic to $\widehat{H}$,
it follows from   Theorem \ref{thm-induced}, that
the twisted  action $(\alpha,\tau): (G_\om,H_\om)\car C_0(\widehat{H})$ is induced from the
twisted action $(\id, \sigma): ((G_H)_\om, H_\om)\car \C$, where $\sigma:H_\om\to \T\cong U(\C)$ is given via evaluation
of $\tau$ at the trivial character $1_H\in \widehat{H}$, i.e., we have $\sigma_{(h,z)}=z$. 
In particular, $\sigma$ is trivial on $H\cong H\times\{1\}\subseteq (G_H)_\om$.
By Theorem \ref{thm-Green}
it follows that $C^*(G,\om)\cong C_0(\widehat{H})\rtimes_{(\alpha,\tau)}G_\om/H_\om$ is Morita equivalent to
$\C\rtimes_{(\id,\sigma)}(G_H)_\om/H_\om\stackrel{(*)}{\cong} C^*(G_H/H, \om')=C^*(G',\om')$, where the isomorphism
$(*)$ follows from the discussion preceding this lemma.  The final assertion now follows from the isomorphism 
$G_\om/(G_H)_\om\cong G/G_H\cong \widehat{H}$ together with \cite{Green-imp}*{Theorem 2.13(i)}.
\end{proof}

Note that by the structure theory for  locally compact abelian groups, every such group is isomorphic 
to some direct product $V\times H$, where $V$ is a vector group of finite dimension $n\geq 0$ and $H$ is a locally compact abelian group 
which contains a compact  open 
 subgroup $K$ (e.g., see \cite{DE}*{Theorem 4.2.1}). Moreover, if $G$ is compactly generated, 
then $H$ is isomorphic to a direct product $Z\times K$, where $Z$ is a finitely generated free abelian group of rank $l\geq 0$
and $K$ is compact. 
We now show that, up to Morita equivalence of the twisted group algebras, we may assume without loss of generality 
that the compact part $K$ is trivial. The  following lemma is analogous to  the first step of the proof of 
\cite{Pogi}*{Theorem 1}  (see \cite{Pogi}*{p.~156}). 

\begin{lemma}\label{lem-K}
Suppose that $\om\in Z^2(G,\T)$ is a totally skew $2$-cocycle on the locally compact abelian group  $G=V\times H$ as above.
Then there exists a discrete abelian group $Z'$ and a totally skew $2$-cocycle $\om'$ on $G'=V\times Z'$ such that
$C^*(G,\om)$ is isomorphic  $C^*(G',\om')\otimes \K(\H)$ for some Hilbert space $\H$.

If, moreover,  $G=V\times Z\times K$ is compactly generated, then $K\cong T\times F$, with $T$ a torus group of dimension
$\dim(T)\leq \rank(Z)$ and $F$ finite, and then 
$Z'$ can be chosen to be free abelian with $\rank(Z')=\rank(Z)-\dim(T)\geq 0$, and there exists a copy of $Z'$ inside $G$ such that
$V\cap Z'=\{0\}$ and 
$\om'=\om_{V\times Z'}$. 
\end{lemma}

%
\begin{proof}
We may assume that $\om$ is normalized. Let $K\subseteq H$ be a compact open subgroup and let $S$ denote the symmetrizer
of $\om_K$. It follows from Remark \ref{rem-twisted-group} that $C^*(K,\om_K)$ is a continuous field over $\widehat{S}$
with simple fibre $C^*(K/S, \om''_{K/S})$. Since all irreducible representations of $C^*(K/S, \om''_{K/S})$ are finite dimensional,
it follows that  $C^*(K/S, \om''_{K/S})$ itself is finite dimensional. Hence $K/S$ is finite, and $S$ is also open in $H$.

By  Lemma \ref{lem-action-omtrivial} we may assume that $\tilde\om:=\om_{G_S}$ is inflated from a normalized cocycle $\om'\in Z^2(G',\T)$, $G':=G_S/S$. In particular, $\om_S=1_{S\times S}$ is trivial. 
As in Lemma \ref{lem-action-omtrivial} we can write
$$C^*(G,\om)\cong C_0(\widehat{S})\rtimes_{(\alpha,\tau)}G_\om/S_\om.$$
and since $C^*(G,\om)$ is simple and  $\widehat{S}$ is discrete (since $S$ is compact), it follows that the action $G_\om\car \widehat{S}$ induced from  $\Res_S^G\circ h_\om:G_\om\to \widehat{S}$ is transitive. 
It follows then from Lemma \ref{lem-action-omtrivial} that $C^*(G,\om)$ is Morita equivalent to $C^*(G',\om')$.
Moreover, since $G/G_S$ is isomorphic to $\widehat{S}$, hence discrete, there clearly exists 
a measurable cross section $G/G_S\to G$.  So we even have $C^*(G,\om)\cong C^*(G',\om')\otimes \K(\ell^2(\widehat{S}))$.
Note that since $\widehat{S}$ is discrete, $G_S$ must contain  $V$ and $G_S/S\cong V\times H_S/S$, where $H_S=G_S\cap H$. 
So the result follows with $Z':=H_S/S$.

Suppose now that $G=V\times Z\times K$ is compactly generated. We apply the above reasoning to $K$.
So let $S$ denote the symmetry group of $\om_K$. Since $K$ and $V$ act 
 trivially on $C^*(K,\om_K)\dach\cong \widehat{S}$, it follows that
$\Res_S^G\circ h_\om: G\to \widehat{S}$ factors through a surjective homomorphism $\psi: Z\to \widehat{S}$. It 
follows that $\widehat{S}$ is isomorphic to a quotient  of $Z$ by the subgroup $Z'=\ker \psi\subseteq Z$.
This quotient is isomorphic to  $\tilde{Z}\times \tilde{F}$ for some finitely generated free abelian group $\tilde{Z}$ with $\rank(\tilde{Z})\leq \rank(Z)$ and some finite group $\tilde{F}$. 
Then $\widehat{S}\cong \tilde{Z}\times \tilde{F}$ implies 
 $S\cong T\times\tilde{F}$ is the product of the torus group $T:=\widehat{\tilde{Z}}$, with 
 $\dim(T)=\rank(\tilde{Z})\leq \rank(Z)$, and the finite group $\tilde{F}\cong \widehat{\tilde{F}}$. 
As observed above, the quotient $K/S$ must be finite as well. Together it  follows that $K$ has a splitting 
$K=T\times F$ for some finite group $F$ containing (an isomorphic copy of) $\tilde{F}$ and $K/S\cong F/\tilde{F}$.

It now follows from  the first part that $C^*(G,\om)$ is isomorphic  to $C^*(G',\om')\otimes \K(\ell^2(\widehat{S}))$ with 
$G'=G_S/S\cong V\times Z'\times F_1$,  with $F_1:=F/\tilde{F}$, and such that the restriction of $\om'$ to $F_1\cong K/S$ 
is totally skew. It follows then from \cite{EchRos}*{Lemma 1} that the splitting $G'\cong V\times Z'\times F_1$ can be chosen 
such that $\om'$ becomes equivalent to the  product cocycle $\om'_{V\times Z'}\otimes \om'_{F_1}$ so that
$$C^*(G',\om')\cong C^*(V\times Z', \om'_{V\times Z'})\otimes C^*(F_1,\om'_{F_1})\cong 
C^*(V\times Z', \om_{V\times Z'})\otimes  M_n(\C)$$
for some  integer $n\geq 0$\footnote{Indeed, one can show $n=\frac{|F_1|}2$.}.  
The group $V\times Z'\subseteq G'$ can even be lifted to a subgroup of $G$ such that $\om'_{V\times Z'}=\om_{V\times Z'}$.
Indeed, if $H'\subseteq G_S$ denotes the inverse image of $V\times Z'\subseteq G_S/K$  under the quotient map, replace 
$V\times Z'$ by any image of a continuous splitting homomorphism $\mathfrak{s}:V\times Z'\to H'\subseteq G$.
Replacing $G'$ with $V\times Z'$ and $\om'$ with $\om_{V\times Z'}$ 
finishes the proof.
\end{proof}

If the group $Z'$ in the above lemma is free abelian (of finite or infinite rank), then we are  able to  reduce the complexity even further.
We need the following observation

\begin{lemma}\label{lem-splitting}
Let $G=V\times Z$ with  $V$ a vector group and $Z$  a free abelian group.
\begin{enumerate}
\item[(a)] If $\varphi:G\to \R$ is a surjective  continuous homomorphism, then 
$H:=\ker\varphi$ has a splitting $H\cong V'\times Z'$ for some vector group $V'\subseteq V$ 
of co-dimension one,  and with $Z'\cong Z$.
\item[(b)] If $\psi:G\to \T$ is a surjective continuous homomorphism, then $H=\ker\psi$ has  a splitting 
$H\cong V'\times \Z v\times Z'$ for some $0\neq v\in V$, some vector group $V'\subseteq V$ of codimension one,  and with $Z'\cong Z$.
\end{enumerate}
\end{lemma}
\begin{proof}
For a proof of (a) let $\{e_i:i\in I\}$ denote a set of generators of $Z$. Since $\varphi|_V:V\to\R$ is surjective, we 
can choose for all $i\in I$ an element $v_i\in V$ such that $\varphi(v_i)=-\varphi(e_i)$. We then obtain 
a new splitting $G\cong V\times Z'$ with generating set  $\{(v_i,e_i): i\in I\}$ of $Z'$ and such that $Z'\subseteq H=\ker\varphi$.
We then get $H=V'\times Z'$ with $V'=\ker\varphi\cap V$ a vector group of co-dimension one on $V$.

Now if $\psi:G\to \T$ is as in (b), then exactly as above we find a new splitting $G\cong V\times Z'$ with $Z'\subseteq H$. 
It follows that there is a splitting $V=V'\times \R v$ for some $0\neq v\in V$ such that $V\cap \ker\psi=V'\times \Z v$ and 
hence $H=V'\times \Z v\times Z'$.
\end{proof}

As in the proof of Poguntke's \cite{Pogi}*{Theorem 1} (see \cite{Pogi}*{pp. ~58--160})  the proof of the following lemma uses induction on the dimension of the vector group $V$. Different from Poguntke's, our result also covers free abelian groups of possibly infinite rank   (see also Example \ref{eq-infinite-torus} below), but the basic ideas are the same.

\begin{lemma}\label{lem-free-abelian}
Suppose $\om$ is a totally skew Borel $2$ cocycle on the group $G=V\times Z$ with $V$ a  nontrivial vector group and 
$Z$ free abelian.
Then there exists a free abelian subgroup $Z'\subseteq G$ with $Z'\cong Z$ such that 
$\om_{Z'}$ is totally skew and $C^*(V\times Z,\om)\cong C^*(G,\om)\cong C^*(Z',\om')\otimes \K(\ell^2(\N))$.
\end{lemma}
\begin{proof}
Suppose that  $\dim(V)=n>0$. We claim that there is a subgroup $G'\subseteq G$ isomorphic to $V'\otimes Z'$ 
with $Z'$ isomorphic to $Z$ and $0\leq \dim(V')<\dim(V)$ such that $C^*(G,\om)\cong C^*(G', \om_{G'})\otimes \K(\ell^2(\N))$.
If this is shown, the result  follows by induction on $\dim(V)$.

In a first step let us assume that
there exist $v,w\in V$ such that $h_\om(v)(w)\neq 1$.
 Let $R=\R w$. Since $H^2(\R,\T)=\{0\}$, we may assume that  $\om_R=1_{R\times R}$ is trivial.
Since every bi-character on the vector group $V$ is of the form $e^{2 \pi i \eta}$ for some bilinear 
 map $\eta:V\times V\to \R$, it follows from $h_\om(v)(w)\neq 1$ that 
 $\Res_R^G\circ h_\om: G\to \widehat{R}$ is surjective.
 By Lemma \ref{lem-splitting} applied to $\varphi:=\Res_R^G\circ h_\om$ it follows that $G$ has a splitting 
 $G=V\times \tilde{Z}$ with $Z'\cong Z$ 
 such that $G_R=\ker\varphi=(V\cap \ker\varphi)\times Z'$. Put $\tilde{V}:=V\cap \ker\varphi\supseteq R$.
 Since there  is  a continuous splitting $G/G_R\to G$, 
Lemma \ref{lem-action-omtrivial} implies that there exists a totally skew 
 cocycle $\om'$ on $G'=G_R/R$ that inflates to $\om_{G_R}$ and such that $C^*(G,\om)\cong C^*(G',\om')\otimes \K(L^2(\widehat{R}))$. If we then choose a splitting  $\tilde{V}:=V'\times R$ we get $G'\cong V'\times Z'\subseteq G$ and 
 $\om'$ identifies with the restriction $\om_{G'}$. Since 
  $\dim(V')=n-2$, the claim follows. 
 
%
%

So assume now that $h_\om(v)(w)=1$ for all $v,w\in V$. This means that the restriction $\om_V$ of $\om$ to $V$ is trivial and  the image of $V$ under the homomorphism 
$$h_\om: G\to \widehat{G}\cong \widehat{V}\times \widehat{Z}$$
projects to $\{0\}\in \widehat{V}\cong V$. Since $h_\om$ is injective (since $\om$ is totally skew), we can find a splitting 
$Z=\Z z\times Z'$ such that the projection of $h_\om(V)$ to $\widehat{\Z z}\cong \T$ is non-trivial, and hence surjective.
Applying Lemma \ref{lem-action-omtrivial}  to $H=\Z z$ and item (b) of Lemma \ref{lem-splitting} 
to $\psi:=\Res_{\Z z}^G\circ h_\om: G\to \widehat{\Z z}\cong \T$, we see
that there exists a totally skew cocycle on $G'=G_{\Z z}/\Z z\cong V'\times \Z v\times Z'$, for  some 
$v\in V$, and some vector subgroup $V'\subseteq V$ of co-dimension one.
 Since $\Z z\times Z'\cong Z$, similar arguments as used in the first case above show that there is a subgroup $G'\subseteq G$
 with $G'\cong V'\times Z$ such that $\om'=\om_{G'}$ is totally skew and $C^*(G,\om)\cong C^*(G',\om')\otimes \K(L^2(\T))$.
 Since $\dim(V')=\dim(V)-1$, this proves the claim.
%
%
\end{proof}

The following example was inspired by some discussion with Christopher Deninger. 

\begin{example}\label{eq-infinite-torus}
To see an interesting  application of Lemma \ref{lem-free-abelian} in the case where $Z$ is free abelian of infinite rank, consider the crossed product $C_0(\R^+)\rtimes_{\tr} \Q^+$, where 
$\R^+$ denotes the multiplicative group of positive real numbers on which the subgroup $\Q^+$
of rational elements acts by multiplication. Since $\Q^+$ is dense in $\R^+$ it follows that $C_0(\R^+)\rtimes_{\tr} \Q^+$ is simple.

Since every $x\in \Q^+$  has a unique expression as 
$x=\prod_{n\in \N} p_n^{l_n}$, where $p_n$ denotes the $n$th prime number and 
$(l_n)_{n\in \N}\in Z:=\sum_{\N}\Z$, it follows that the
natural logarithm 
$\ln:\R^+\to \R$ transforms the action $\tr:\Q^+\car \R^+$ to the action 
$$\tilde{\tr}:Z \car \R; \quad ((l_n)_n,y)\mapsto y+\sum_{n\in\N} \ln(p_n)l_n.$$
Now let $G=\R\times Z$ and define $\om\in Z^2(G,\T)$ by
$$\om\big((v, (l_n)_n), (w, (k_n)_n)\big)=\prod_{n\in \N} \exp\big(2\pi i \ln(p_n) l_n w\big).$$
Notice that this is the restriction of the Heisenberg cocycle 
$$\eta:\R\times\R\to \T; \eta((x,y), (x', y'))=
\exp(2\pi i yx')$$
 to $\R\times Z$ via the embedding $Z\into\R; (l_n)_n\mapsto \sum_n \ln(p_n)l_n$.
It then easy to check that 
$$C_0(\R^+)\rtimes_{\tr} \Q^+\cong C_0(\R)\rtimes_{\tilde\tr} Z\cong C^*(\R\times Z, \om).$$
By simplicity, it follows that $\om$ is totally skew, and then Lemma \ref{lem-free-abelian} implies 
that there is totally skew $2$-cocycle $\om'\in Z^2(Z,\T)$ such that 
$C_0(\R^+)\rtimes_{\tr} \Q^+$ is isomorphic to  $C^*(Z,\om')\otimes \K(\ell^2(\N))$.
\end{example}

The following  version of Poguntke's \cite{Pogi}*{Theorem 1} now follows from a straightforward combination of 
Lemma \ref{lem-K} and  Lemma \ref{lem-free-abelian}. 

 \begin{theorem}\label{thm-totallyskew}
 Suppose that $G\cong V\times Z\times K$ is a compactly generated locally compact abelian group and that 
 $\om\in Z^2(G,\T)$ is a totally skew Borel $2$-cocycle on $G$. Then $K\cong T\times F$ for some torus group $T$ with $\dim(T)\leq \rank{Z}$ and $F$ finite, and there exists a copy  $Z'\subseteq G$ of $Z$ such that the restriction 
 $\om':=\om_{Z'}$  is totally skew and such that $C^*(G,\om)\cong C^*(Z',\om')\otimes \K(\H)$. 
Moreover, $\H\cong \ell^2(\N)$ whenever $V\times T$ is nontrivial. Otherwise $\H$ is finite dimensional.
\end{theorem}

\begin{remark}\label{rem-Pogi1}
In \cite{Pogi}*{Theorem 1} Poguntke shows that for every compactly generated two-step nilpotent group $G$ every
simple quotient of $C^*(G)$ is Morita equivalent to some (possibly zero-dimensional) non-commutative torus
 $C^*(\Z^l,\om)$.
This result can easily be deduced from Theorem \ref{thm-totallyskew} above by observing that in this case
$C^*(G)$ is an algebra of sections of a continuous field of $C^*$-algebras over $\widehat{Z}$, where 
$Z$ denotes the center of $G$, with fibre over $\chi\in \widehat{Z}$ isomorphic to a twisted group algebra
$C^*(G/Z, \om_\chi)$, for some cocycle $\om_\chi\in Z^2(G/Z,\T)$ depending on $\chi$.
Note that $G/Z$ is compactly generated abelian if $G$ is compactly generated two-step nilpotent.

Conversely, Theorem \ref{thm-totallyskew}  is of course also an easy consequence of Poguntke's \cite{Pogi}*{Theorem 1}
applied  to $C^*(G_\om)$, if $G_\om$ denotes the central extension of $G$ by $\T$ corresponding to $\om$.
%
%

\end{remark}

\begin{remark}\label{rem-compactly-generated}
By  Lemma \ref{lem-K} we see that for abelian $G$ and  a totally skew $2$-cocycle $\om \in Z^2(G,\T)$  the twisted group algebra $C^*(G,\om)$ is always Morita equivalent to
a twisted group algebra $C^*(V\times D,\om')$ for some vector group $V$ and some discrete group $D$. One might wonder 
whether one could use arguments similar to the ones given in the proof of Lemma  \ref{lem-free-abelian}, which allow 
a reduction to the case of discrete groups. A good test case would be the twisted group algebra 
$C^*(\R\times\mathbb{Q},\om)$ with totally skew cocycle $\om((x,y), (x', y'))=e^{2\pi i xy'}$, for which an analogue of the induction 
step as in the proof of  Lemma  \ref{lem-free-abelian} certainly fails.

\end{remark}

\section{Conclusions}\label{sec-Poguntke}

We first combine Theorems \ref{thm-main}  together with Theorem \ref{thm-totallyskew}, to get

\begin{corollary}\label{cor-torus}
Let $\gamma:L\car A$ and $\alpha=\gamma \circ \varphi:G\car A$ be as in Theorem \ref{thm-main} such that $G$ is abelian. 
Then the following are true:
\begin{enumerate}
\item All point-sets $\{P\}\subseteq \Prim(A\rtimes_{\alpha}G)$ are locally closed.
\item Every simple subquotient $B_P$ of $B:=A\rtimes_\alpha G$ is Morita equivalent  to some twisted group algebra $C^*(Z_P,\om_P)$ of some abelian locally compact group $Z_P$ (which could be trivial) and some totally skew Borel $2$-cocycle $\om_P$ on $Z_P$. 
\item  If $G$ and all stabilizers $H_\rho$ as in Theorem \ref{thm-main} are 
compactly generated, then all $Z_P$ can be chosen to be either trivial or finitely generated  free abelian. Hence $C^*(Z_P,\om_P)$ is  either $\C$ or  a simple 
non-commutative torus of dimension $\geq 2$.
\end{enumerate}
\end{corollary}

The following corollary is Poguntke's result for 
connected groups

\begin{corollary}\label{cor-connected}
Suppose that $G$ is a connected locally compact group. Then every one-point set $\{P\}\subseteq \Prim(C^*(G))$ is locally closed 
and the corresponding simple subquotient $C^*(G)_P$ is either isomorphic to an algebra of compact operators $\K(\H)$, or it  is Morita equivalent to a simple non-commutative torus of dimension $\geq 2$.
\end{corollary}
\begin{proof}

In case that $G$ is a connected and simply connected Lie group, the result follows from Corollary \ref{cor-simplyconnected} 
and Theorem \ref{thm-totallyskew}. For general connected groups, we follow Poguntke's arguments in order to reduce to this case.

By the structure theory for connected groups (e.g., see \cite{MontZip}) we know that there exists a compact normal subgroup 
$K$ of $G$ such that $\dot{G}:=G/K$ is a Lie group. We can then write $C^*(G)$ as a twisted crossed product 
$C^*(K)\rtimes_{(\alpha,\tau)}G/K$, and then, up to Morita equivalence, as an ordinary  crossed product 
$B\rtimes_\beta \dot{G}$, with $B$ Morita equivalent to $C^*(K)$. 
Since $K$ is compact, $C^*(K)$ is type I with discrete dual $\widehat{K}=\widehat{C^*(K)}$, and hence $B$ is type I with 
discrete dual $\widehat{B}$. Since $\dot{G}$ is connected, the action of $\dot{G}$ on $\widehat{B}$ is trivial, 
and we therefore obtain a decomposition
$$B\rtimes_\beta \dot{G}\cong \bigoplus_{\rho\in \widehat{B}}\K(\H_\rho)\rtimes_{\beta^\rho}\dot{G}\cong \bigoplus_{\rho\in \widehat{B}} \K(\H_\rho)\otimes C^*(\dot{G},\om_\rho)$$
with $\om_\rho\in Z^2(\dot{G}, \T)$ a suitable $2$-cocycle for every $\rho\in \widehat{B}$. 
Thus, the conclusion of Theorem \ref{thm-Pogi} will follow if we can show that the same conclusion holds for 
all twisted group algebras $C^*(\dot{G},\om)$ for $\om\in Z^2(\dot{G}, \T)$. But we observed in 
Remark \ref{rem-quotient} that $C^*(\dot{G},\om)$ is a quotient (indeed a direct  summand) of $C^*(\dot{G}_\om)$, 
where $\dot{G}_\om$ is the central extension of $\dot{G}$ by $\T$ corresponding to $\om$. Since this is a connected Lie group, 
the conclusion of Theorem \ref{thm-Pogi} will follow for all connected groups, if we can show it for 
connected Lie groups.

 Next, if $G$ is a connectetd Lie group and $\tilde{G}$ is the universal covering group of $G$, then  $G=\tilde{G}/Z$ for some discrete central subgroup $Z$ of $G$.
 It follows that $C^*(G)$ is  a quotient 
of $C^*(\tilde{G})$. Thus every simply subquotient of $C^*(G)$ is also a simple subquotient of $C^*(\tilde{G})$ and the result follows.
%
%
%
\end{proof}

The proof of Theorem \ref{thm-main-intro} now follows from a combination of Corollary \ref{cor-connected} with Theorem \ref{thm-totallyskew} (see also Remark \ref{rem-Pogi1}). To illustrate that our results go strictly further than just the case of connected groups, we consider the following

\begin{example}\label{ex-Mautner}
Recall  that for any irrational number $\theta$, the  {\em Mautner group} $M_\theta:=\C^2\rtimes_\theta \R$ is the semi-direct product of the additive 
group $\C^2$ by the real line $\R$ 
acting on $\C^2$ via the automorphisms $t\cdot(z,w):=(e^{2\pi i t}z, e^{ 2\pi i\theta t}w)$.
It is a well-known fact that $M_\theta$ is a connected Lie group with minimal  dimension such that $C^*(M_\theta)$ is not type I. 

Instead of regarding this example via the general result on connected groups, we shall  apply Corollary \ref{cor-compact} instead.
Since $M_\theta=\C^2\rtimes_\theta \R$ we have 
$$C^*(M_\theta)\cong C^*(\C^2)\rtimes_{\alpha^\theta} \R\cong C_0(\widehat{\C^2})\rtimes_{\alpha^\theta_*}\R,$$
where $\alpha^\theta$ denotes the decomposition action and $\alpha^\theta_*$ is the action on $\widehat{\C^2}$ given by the formula $(\alpha^\theta_*)_t(\chi)(z,w)=\chi\big((-t)\cdot(z,w))$ for $\chi\in \widehat{\C^2}$.
Using the canonical identifications  $\widehat{\C^2}\cong\widehat{\R^4}\cong \R^4\cong\C^2$ a little exercise shows that 
$$C^*(M_\theta)\cong  C_0(\widehat{\C^2})\rtimes_{{\alpha^\theta_*}}\R\cong C_0(\C^2)\rtimes_{\alpha}\R$$
with 
\begin{equation}\label{eq-actionbeta}
\alpha_t(f)(z,w)=f(t\cdot(z,w))=f(e^{2\pi it}z, e^{2\pi i \theta t}w), \quad f\in C_0(\C^2).
\end{equation}
Consider now the compact action $\gamma:\T^2\car C_0(\C^2)$ induced from the action $\T^2\car\C^2$ given by 
$(u,v)\cdot (z,w)=(\bar{u}v, \bar{z}w)$,  for $({u},v)\in \T^2$, and $(z,w)\in \C^2$.  Then $\alpha=\gamma\circ \varphi$ with 
$\varphi:\R\to \T^2, \varphi(t)=(e^{2\pi i t},  e^{2\pi i \theta t})$. 
The orbit space $\C^2/\T^2$ is homeomorphic to 
$[0,\infty)^2$ with orbits 
\begin{equation}\label{eq-orbit}
\mathcal O_{x,y}=\{(z,w)\in \C^2: |z|=x, |w|=y\},\quad (x,y)\in [0,1)^2.
\end{equation}
If $x\neq 0\neq y$, the orbit $\mathcal O_{x,y}$ is  homeomorphic to $\T^2$ with actions as above. The other orbits are  
homeomorphic to $\{0\}\times \T$  (if $x=0\neq y$), to $\T\times\{0\}$ (if $x\neq 0=y$) and  we have $\mathcal O_{(0,0)}=\{(0,0)\}$.

By Theorem \ref{thm-main} we see that $C^*(M_{\theta})\cong C_0(\C^2)\rtimes_\alpha \R$ has a decomposition by subquotients 
$C(\mathcal{O}_{x,y})\rtimes_\alpha \R$ over $[0,1)^2$,  each of them isomorphic to twisted group algebras.
If $(x,y)=(0,0)$, we obtain  $C^*(\R)\cong C_0(\R)$. If $x=0\neq y$ and $x\neq 0=y$, we get copies of 
$C(\T)\rtimes \R\sim_M C^*(\Z)\cong C(\T)$.
So all these fibres are type I. So let us study the case $x\neq 0\neq y$, in which case we have 
$C(\mathcal O_{x,y})\rtimes_\alpha\R\cong C(\T^2)\rtimes_{\alpha}\R$.
Identifying $\T^2$ with $(\R/\Z)\times \T$, we have the  $\R$-equivariant  projection $p: \T^2\to \R/\Z$ to the first factor, from which it follows by 
Theorem \ref{thm-induced}  together with Green's Theorem \ref{thm-Green} that 
$$C(\T^2)\rtimes_\alpha\R\cong \Ind_{\Z}^{\R}C(\T)\rtimes_{\Ind \alpha^\theta}\R\sim_M C(\T)\rtimes_{\alpha^\theta}\Z$$
where $\alpha^\theta:\Z\car C(\T)$ denotes rotation by the irrational angle $\theta$. 
This is the simple non-commutative 2-torus $A_\theta$.
\medskip

There are many variations of this example: let $G$ be any abelian locally compact group and let $\chi,\mu \in \widehat{G}$ be any two characters of $G$. We can then form the semidirect product $M_{\chi,\mu}:=\C^2\rtimes^{\chi,\mu} G$ 
with respect to the action $g\cdot(z,w)=(\chi(g)z, \mu(g)w)$, $g\in G$, $(z,w)\in \C^2$. Then Corollary \ref{cor-compact} 
and Theorem \ref{thm-main} will apply to obtain similar decompositions of $C^*(M_{\chi,\mu})$ by 
subquotients $C(\mathcal{O}_{x,y})\rtimes_\alpha G$ as above, which  for $x\neq 0\neq y$ turn out to be Morita equivalent  
to twisted group algebras $C^*(G\times \Z^2, \om)$ for the cocycle 
$$\om\big((g_1, n_1, n_2), (g_2, n_2, m_2)\big) = \chi(g_1)^{n_2}\mu(g_1)^{m_2}.$$
If $G$ is compactly generated, then so is $G\times\Z^2$ and hence all non-elementary simple quotients\footnote{It is not difficult to see that all primitive ideals of $C^*(M_{\chi,\mu})$ are maximal.} of $C^*(M_{\chi,\mu})$ are Morita equivalent to non-commutative tori.
We omit further details. Note that the $M_{\chi,\mu}$ are  covered by Poguntke's \cite{Pogi}*{Theorem 2} if and only if $G$ is connected.

\end{example}

\section{The case of almost connected groups}\label{sec-almost connected}

Recall that a locally compact group is called {\em almost connected} if the quotient $G/G_0$ of $G$ by its  connected component $G_0$ is compact.
The aim of this section is to study simple subquotients of almost connected groups. As for connected groups, it is well known that every almost connected group 
 has a compact normal subgroup $N$ such that $G/N$ is a (almost connected)  Lie group.
Precisely the same arguments as used in the proof of Corollary \ref{cor-connected} then show that, in order to study the structure of simple subquotients of 
group $C^*$-algebras of  almost connected groups, it suffices  to study the case of Lie groups, where, in particular, the quotient group $G/G_0$ is finite.

We may then write $C^*(G)$ as a twisted crossed product $C^*(G_0)\rtimes_{(\alpha,\tau)}G/G_0$, and by passing to a Morita equivalent ordinary action 
$\beta:F\car B$ for the quotient group $F:=G/G_0$, we are left to study (up to Morita equivalence) crossed products $B\rtimes_\beta F$ by finite group actions,
where $B$ is Morita equivalent to the group algebra $C^*(G_0)$.
By the results of the previous section, we then know that 
\begin{itemize}
\item all points $\{P\}\in \Prim(B)$ are locally closed, and
\item the simple subquotients $B_P$ of $B$ corresponding to some $P\in \Prim(B)$ are either Morita equivalent 
to $\C$ or they are Morita equivalent to some simple non-commutative  $C^*(Z,\om)$ of dimension $n\geq 2$.
\end{itemize}

We are now going to show the following:

\begin{theorem}\label{thm-crossed}
Suppose that $\beta:F\car B$ is an action of a finite group $F$ on the separable $C^*$-algebra $B$ such that 
every one-point set $\{P\}$ in $\Prim(B)$ is locally closed. For every $P\in \Prim(B)$ let $\beta^P:F_P\car B_P$ denote the 
action of the stabilizer $F_P$ of $P$ on the simple subquotient $B_P$ induced by $\beta$.
Then the following are true:
\begin{enumerate}
\item All orbits $F(P)$ in $\Prim(B)$ are locally closed. 
\item Every one-point set $\{Q\}\in \Prim(B\rtimes_\beta F)$ is locally closed.
\item The simple subquotient $(B\rtimes_\beta F)_Q$  corresponding to  \mbox{$Q\in\Prim(B\rtimes_\beta F)$} is Morita equivalent  to 
a simple subquotient of some crossed product $B_P\rtimes_{\beta^P}F_P$.
\end{enumerate}
\end{theorem}

The result will follow almost immediately from the Mackey-Rieffel-Green machine, once we have the following

\begin{lemma}\label{lem-orbit}
Suppose that $F$ is a finite group acting by homeomorphisms on the topological T$_0$-space $X$. 
Then, for every $x\in X$,  the orbit $F(x)=\{gx:g\in F\}$ is discrete in the relative topology. 
If, moreover,  $\{x\}$ is locally closed in $X$, then so is $F(x)$.
\end{lemma}
\begin{proof}
For the  first assertion, we first show that for every $x\in X$  the orbit $F(x)=\{gx: g\in F\}$ is a T$_1$-space in the relative topology.
Indeed, if this fails to be true, there exists some $x\in X$ and some $g\in F$ with $gx\neq x$, but $gx\in \overline{\{x\}}$.
We claim that this implies that $x\in \overline{\{gx\}}$ as well, and hence that $gx=x$, since $X$ is a T$_0$-space.
Indeed, since $gx\in \overline{\{x\}}$, it follows that $g^nx\in \overline{\{gx\}}$ for all $n\geq 1$. This is clear for 
$n=1$ and if it is true for $n$, we also get 
$$g^{n+1}x=g^n(gx)\subseteq g^{n}\overline{\{x\}}= \overline{\{g^nx\}}\subseteq\overline{\{gx\}}.$$
Since $F$ is finite, there exists some $n\in \N$ with $g^nx=x$, hence $x\in \overline{\{gx\}}$. 
It follows from this that every subset of $F(x)=\{gx:g\in F\}$ is relatively closed, and hence that $F(x)$ is discrete with respect to the relative topology.
The second assertion follows from the fact that finite unions of locally closed sets are locally closed.
\end{proof}

\begin{proof}[Proof of Theorem \ref{thm-crossed}]
The first assertion is a direct consequence of Lemma \ref{lem-orbit}. 
It follows then from an application of Proposition \ref{orbit-method-separable} that for every $Q\in \Prim(B\rtimes_\beta F)$, there exists a unique orbit 
$F(P)\subseteq \Prim(B)$ such that $Q$ is an element of the locally closed subset $\Prim(B_{F(P)}\rtimes_{\beta}F)$.
Since, again by Lemma \ref{lem-orbit}, the set $F(P)$ is discrete in the relative topology, it follows that $F(P)$ is canonically homeomorphic to
$F/F_P$, and hence $(B_{F(P)}, F, \beta)$ is isomorphic to $(\Ind_{F_P}^FB_P, F, \Ind\beta^P)$. By Green's  Theorem \ref{thm-Green} it follows 
that $B_{F(P)}\rtimes_{\beta}F$ is Morita equivalent to $B_P\rtimes_{\beta^P}F_P$.
Since $B_P$ is simple and $F_P$ is finite, it  follows then from \cite{Ech:T1}*{Theorem 3.4} (or from a crossed-product version of \cite{Pogi:finite}*{Satz 2}) that every primitive 
ideal in $B_P\rtimes_{\beta^P}F_P$ is maximal. By Morita equivalence,  the same holds then for $B_{F(P)}\rtimes_{\beta}F$.
Hence $\{Q\}$ is closed in the locally closed subset $\Prim(B_{F(P)}\rtimes_{\beta}F)$ of $\Prim(B\rtimes_\beta F)$. 
By Proposition \ref{prop-loc-closed-points}, it follows that $\{Q\}$ is locally closed in $\Prim(B\rtimes_\beta F)$.
This finishes the proof.
\end{proof}

\begin{corollary}\label{cor-almost-connected}
Let $G$ be an almost connected locally compact group. Then every one-point set $\{Q\}$ in $\Prim(C^*(G))$ is locally closed.
Moreover, the simple subquotient of $C^*(G)$ corresponding to $Q$ is either isomorphic to $\K(\H)$ for some Hilbert space $\H$
or it is Morita equivalent to some simple quotient of a crossed product of the form $(C^*(Z,\om)\otimes\K(\H))\rtimes_\beta F$
for some simple non-commutative $n$-torus $C^*(Z,\om)$ and some finite group $F$.
\end{corollary}
\begin{proof} This follows from Theorem \ref{thm-crossed} together with the discussion preceding it.
\end{proof}

It would  therefore be very interesting to study the structure of crossed products $(C^*(Z,\om)\otimes\K(\H))\rtimes_\beta F$
with $C^*(Z,\om)$ a simple non-commutative $n$-torus of dimension $n\geq 2$ and $F$ a finite group.
Some work has been done already for the case of actions $\beta:F\car C^*(Z,\om)$  (e.g., see
\cites{ELPW, He, Jeong, Chakraborty}), but the general situation is still very unclear.
We close with the following example, which might illustrate what one could expect:

\begin{example}\label{ex-Mautner-almost}
Recall  from Example \ref{ex-Mautner} the  Mautner group $M_\theta:=\C^2\rtimes_\theta \R$.
In what follows, let us write the elements of $M_\theta$ as triples $(z,w, t)$ with $(z,w)\in \C^2$, $t\in \R$. 
There is an action of the two-element group  $\Z_2:=\Z/2\Z$ on $M_\theta$ given by
$$\epsilon\cdot(z,w,t)=(\bar{z}, \bar{w}, -t)$$
for the generator $\epsilon \in \Z_2$. 
We then obtain the almost connected group $M_\theta\rtimes \Z_2$ which can also be realized as a semi-direct product 
$\C^2\rtimes(\R\rtimes \Z_2)$.
As in Example \ref{ex-Mautner} we get
$$C^*(M_\theta\rtimes \Z_2)\cong C^*(\C^2\rtimes(\R\rtimes \Z_2))\cong C^*(\C^2)\rtimes_{\alpha_*}(\R\rtimes\Z_2)\cong C_0(\C^2)\rtimes_\alpha(\R\rtimes\Z_2)$$
with action $\alpha:\R\rtimes\Z_2\car C_0(\C^2)$ given for $t\in \R$ as in Example \ref{ex-Mautner} and on the generator $\epsilon\in \Z_2$ by 
$\epsilon\cdot (z,w)=(\bar{z}, \bar{w})$. Now consider the $\T^2$-orbits $\mathcal O_{x,y}$ in $\C^2$ as in (\ref{eq-orbit}), which 
are all invariant under the action of $\R\rtimes\Z_2$. We therefore may decompose 
$\Prim(C^*(M_\theta\rtimes \Z_2))$ into the pairwise disjoint union of the spaces $\Prim(C(\mathcal O_{x,y})\rtimes_\alpha \R\rtimes\Z_2)$.
If $x=0$ or $y=0$, then these crossed products are easily seen to be type I, so let us restrict to the case $x\neq 0\neq y$.
As in Example \ref{ex-Mautner} we can then  identify $\mathcal O_{x,y}$ with $\T^2$. The action of $\R\rtimes \Z_2$ on the first factor $\T$ of $\T^2$
has stabilizer at $1\in \T$ given by $\Z\rtimes \Z_2$, so we get a canonical $\R\rtimes\Z_2$-equivariant homeomorphism 
$\T\cong (\R\times\Z_2)/(\Z\rtimes\Z_2)$. We therefore obtain an $\R\rtimes\Z_2$-equivariant first-factor  projection $p:\T^2\to \T\cong  (\R\times\Z_2)/(\Z\rtimes\Z_2)$.
Hence by Theorem \ref{thm-induced} and Green's Theorem \ref{thm-Green} we get
\begin{align*}
C(\T^2)\rtimes_\alpha(\R\times\Z_2)&\cong \Ind_{\Z\rtimes\Z_2}^{\R\rtimes\Z_2}C(\T)\rtimes_{\Ind \alpha^\theta}(\R\rtimes\Z_2)\\
&\sim_M C(\T)\rtimes_{\alpha^\theta}(\Z\rtimes\Z_2)\\
&\cong \big(C(\T)\rtimes_{\alpha^\theta}\Z\big)\rtimes\Z_2\\
&\cong A_\theta\rtimes \Z_2
\end{align*}
where the action $\Z_2\car A_{\theta}$ on the irrational rotation algebra $A_{\theta}$ 
is given on the  standard unitary generators $u, v$ of $A_\theta$ by $u\mapsto u^*, v\mapsto v^*$.
Thus we see that the crossed products $C(\mathcal O_{x,y})\rtimes_\alpha (\R\rtimes\Z_2)$ for $x\neq 0\neq y$ are all stably isomorphic to
(indeed,  they are stabilizations of) the crossed products $A_\theta\rtimes \Z_2$ by the flip automorphism. These have been shown to 
be non type I AF-algebras by Bratteli and Kishimoto in \cite{BK} (see also \cite{ELPW}).
\end{example}

\end{document}